\RequirePackage{fix-cm}
\documentclass[10pt,a4paper]{article}
\usepackage[margin=20mm]{geometry}
\usepackage{amsmath,amssymb,amsthm,booktabs,tikz,microtype}
\usepackage[inline]{enumitem}
\setlist[enumerate]{align=left}
\allowdisplaybreaks[4]
\usepackage{xcolor}
\usepackage[colorlinks=true,linkcolor=blue,citecolor=blue,urlcolor=blue]{hyperref}
\usetikzlibrary{positioning,fit}
\tikzset{vertex/.style={circle,fill,draw=none,inner sep=0pt,outer sep=0pt,minimum size=3.2pt,transform shape=false}}
\numberwithin{equation}{section}
\newtheorem{theorem}{Theorem}[section]
\newtheorem{corollary}[theorem]{Corollary}
\newtheorem{lemma}[theorem]{Lemma}
\newtheorem{claim}[theorem]{Claim}
\newtheorem{observation}[theorem]{Observation}
\theoremstyle{definition}
\newtheorem{construction}[theorem]{Construction}
\newtheorem{remark}[theorem]{Remark}
\newcommand{\ov}[1]{\overline{#1}}
\newcommand{\pp}[1]{\left[#1\right]_{+}}

\title{The maximum number of triangles in regular graphs}
\author{\fontsize{11pt}{13.2pt}\selectfont Ya-Hong Chen$^{a}$,\quad Shi-Cai Gong$^{b}$,\quad Changhong Lu$^{c}$, \quad Zhuojie Tu$^{c}$,\quad Long-Tu Yuan$^{c}$\\[6pt]
\small $^{a}$Department of Mathematics, Lishui University,  Lishui, 323000, PR China.\\[3pt]
\small $^{b}$School of Science, Zhejiang University of Science and Technology, Hangzhou, 310023, PR China.\\[3pt]
\small $^{c}$School of Mathematical Sciences, Key Laboratory of MEA (Ministry of Education),\\ 
\small Shanghai Key Laboratory of PMMP,
 and Nantong Institute for Applied Mathematics and Artificial Intelligence,\\
\small East China Normal University, Shanghai, 200241, China.}
\date{}
\begin{document}
\maketitle
\vspace{-1.5em}
\setlength{\parskip}{0pt}
\begingroup\renewcommand{\thefootnote}{}\footnotetext{E-mail addresses: \texttt{yhchen0611@163.com} (Ya-Hong Chen); \texttt{scgong@zafu.edu.cn} (Shi-Cai Gong); \texttt{chlu@math.ecnu.edu.cn} (Changhong Lu); \texttt{tu$_{-}$zhuojie@163.com} (Zhuojie Tu); \texttt{ltyuan@math.ecnu.edu.cn} (Long-Tu Yuan).}\endgroup
\begin{abstract}
We determine the maximum number of triangles in $n$-vertex $k$-regular graphs and characterize all extremal graphs for all admissible orders when $k$ is odd and for an explicit range of orders when $k$ is even. For even $k$, we also give explicit lower and upper bounds for the remaining orders. As a consequence, we obtain a lower bound on the number of independent sets of size three in $k$-regular graphs and characterize all equality cases. Our exact results confirm the conjecture of Petingi and Rodr\'iguez and the triangle case of a conjecture of Cambie, de Joannis de Verclos and Kang for all odd $k$ and the specified orders when $k$ is even. Moreover, for every even $k\ge18$, we construct counterexamples to both conjectures and obtain a recursive characterization of all extremal graphs of order at least $Ck$, where $C$ is an absolute constant.
\end{abstract}

\section{Introduction}
Maximizing the number of copies of a fixed graph in a graph on $n$ vertices  with maximum degree at most $k$ is a special case of the generalized Tur\'an problems studied by Alon and Shikhelman~\cite{AS}. Indeed, the degree constraint is equivalent to forbidding the star $K_{1,k+1}$ as a subgraph. Clique counting under a maximum degree constraint has received particular attention. Denote by $K_t$ the \emph{complete graph} on $t$ vertices. Gan, Loh and Sudakov~\cite{GLS} conjectured that, for each fixed $t\ge3$, disjoint unions of complete graphs maximize the number of copies of $K_t$ under a maximum degree constraint, and reduced the conjecture to the triangle case. Cutler and Radcliffe~\cite{CRfixed} verified the conjecture for maximum degree at most 6, and Chase~\cite{Chase} settled the triangle case of the conjecture. Specifically, if $n=q(k+1)+r$ for integers $q,r$ with $q\ge0$ and $0\le r<k+1$, then every graph on $n$ vertices with maximum degree at most $k$ has at most $q\binom{k+1}{3}+\binom r3$ triangles, and this bound is attained by $q$ disjoint copies of $K_{k+1}$, together with $K_r$ when $r>0$. Note that this extremal graph is not $k$-regular when $0<r<k+1$. This leads to the corresponding extremal counting problem for \emph{$k$-regular} graphs, in which every vertex has degree exactly $k$. Extremal counting problems for other types of subgraphs in regular graphs have also been studied~\cite{DJPR,Harangi,McKay,Schrijver} and here we focus on the maximum number of triangles in regular graphs of prescribed order and degree. 

The additional requirement that the graph be $k$-regular can strictly reduce the maximum number of triangles. For example, for every positive even $k$, results of Cambie, de Joannis de Verclos and Kang~\cite{CDK} imply that requiring every vertex to have degree $k$ reduces the maximum number of triangles on $2k+1$ vertices by exactly $k^2/4$ compared with allowing maximum degree at most $k$. A special case of their~\cite{CDK}  conjecture states that, when $n=q(k+1)+r$ with $0<r<k+1$, every $k$-regular graph on $n$ vertices with the maximum number of triangles has at least $q-1$ components isomorphic to $K_{k+1}$. In this paper, we call this special case the \textbf{CDK Conjecture}. The intuition behind this conjecture is also reflected in an earlier result of van der Hoorn, Lippner and Mossel~\cite{HLM}, who showed that, for fixed $k\ge3$ and fixed $0<c<1$, almost all labelled $k$-regular graphs on $n$ vertices with at least $cn\binom k2/3$ triangles have all but $o(n)$ of their triangles in components isomorphic to $K_{k+1}$.

In this paper, we determine the maximum number of triangles and characterize all extremal graphs for every odd $k$ and for certain orders when $k$ is even, as stated in Theorem~\ref{mainodd}. To state this result, we first introduce some notation. Let $\tau(G)$ be the number of copies of $K_3$ in a graph $G$. For integers $n>k\ge0$ with $nk$ even, define $\tau_k(n)$ as the maximum of $\tau(G)$ over all $k$-regular graphs $G$ on $n$ vertices. A $k$-regular graph $G$ on $n$ vertices is called \emph{extremal} if $\tau(G)=\tau_k(n)$. We write $\ov H$ for the \emph{complement} of $H$ and $G\cup H$ for the \emph{disjoint union} of $G$ and $H$. In particular, $mH$ denotes the disjoint union of $m$ copies of $H$.

In the following theorems and corollaries, we consider $k$-regular graphs on $n$ vertices, where $k\ge1$, $n\ge k+1$ and $nk$ is even. Write $n=q(k+1)+r$, where $q\ge1$ and $0\le r<k+1$. Let $$f(n,k)=\frac{nk(k-1)-r(k+1-r)(k+1+r)}{6}.$$

\begin{theorem}\label{mainodd}
If $k$ is odd, or $r$ is odd, or $3r\le2(k+1)$, then $\tau_k(n)=f(n,k).$
Moreover, the extremal graphs are exactly the following.
\begin{enumerate}[label=(\roman*),leftmargin=2em]
\item [$(i).$] For $r=0$, the graph $qK_{k+1}$.
\item [$(ii).$] For $r>0$, the graphs $(q-1)K_{k+1}\cup\ov H$, where $H$ is any triangle-free $r$-regular graph on $k+1+r$ vertices.
\end{enumerate}
\end{theorem}

For the remaining orders, we also give the following explicit lower and upper bounds.

\begingroup
\begin{theorem}\label{mainbounds}
If $k$ and $r$ are even and $3r>2(k+1)$, then
\[
\begin{aligned}
f(n,k)-\frac{r(2r-k-2)}4
\le\tau_k(n)\leq f(n,k) - 1.
\end{aligned}
\]
\end{theorem}
\endgroup


For even $k\ge18$, the orders $3k$ and $2k-1$ both have the even remainder $r=k-2$ modulo $k+1$, with $3r>2(k+1)$, and hence lie outside the range of Theorem~\ref{mainodd}. The following theorem shows that, at order $3k$, no extremal graph contains a component isomorphic to $K_{k+1}$.
\begin{theorem}\label{maineven}
For every even integer $k\ge18$,
\[
\tau_k(3k)\ge f(3k,k)-\frac{3k}{2}+1
\quad\text{and}\quad
\tau_k(2k-1)\le f(2k-1,k)-\frac{3k}{2}.
\]
In particular, $\tau_k(3k)>\tau_k(2k-1)+\binom{k+1}{3}$.
\end{theorem}

\begingroup
Nevertheless, for every even $k\ge18$ and all sufficiently large $n$, the following theorem gives a recurrence for the maximum number of triangles in $k$-regular graphs and a corresponding recursive characterization of the extremal graphs.
\begin{theorem}\label{mainperiodicity}
For every even integer $k\ge18$ and every integer $n\ge Ck$, we have
\[
\tau_k(n+k+1)=\tau_k(n)+\binom{k+1}{3},
\]
where $C$ is an absolute constant. Moreover, the extremal $k$-regular graphs on $n+k+1$ vertices are precisely those obtained by adjoining a disjoint $K_{k+1}$ to an extremal $k$-regular graph on $n$ vertices.
\end{theorem}
\endgroup

\begingroup
Theorems~\ref{mainodd} and~\ref{maineven} yield the following corollary concerning the CDK Conjecture.
\begin{corollary}\label{cdk_consequences}
The \textbf{CDK Conjecture} holds whenever $k$ is odd, or $r$ is odd, or $3r\le2(k+1)$. However, for every even $k\ge18$ and every integer $m\ge0$, there is a counterexample of order $3k+m(k+1)$.
\end{corollary}
\begin{proof}
The first assertion follows from Theorem~\ref{mainodd}. For the second, suppose that an extremal graph on $n=3k+m(k+1)=(m+2)(k+1)+(k-2)$ vertices contains $m+1$ components isomorphic to $K_{k+1}$, and let $H$ be the graph obtained by deleting these $m+1$ components. Since $|V(H)|=2k-1$, Theorem~\ref{maineven} shows that replacing $H$ and one of these components with an extremal graph on $3k$ vertices increases the number of triangles, a contradiction.
\end{proof}
\endgroup

Furthermore, for fixed order and degree, triangle maximization in regular graphs is equivalent to minimizing the number of induced paths on three vertices. Denote by $P_n$ the \emph{path} on $n$ vertices. Let $\nu(G)$ be the number of induced copies of $P_3$ in a graph $G$, and for integers $n>k\ge0$ with $nk$ even, define $\nu_k(n)$ as the minimum of $\nu(G)$ over all $k$-regular graphs $G$ on $n$ vertices. To see this equivalence, let $G$ be a $k$-regular graph on $n$ vertices. When we count pairs of neighbors of each vertex, each induced copy of $P_3$ is counted once at its middle vertex, while each triangle is counted once at each of its three vertices. Thus
\begin{equation}\label{intro_path_identity}
\nu(G)+3\tau(G)=n\binom{k}{2},
\end{equation}
and hence
\begin{equation}\label{intro_path_minimum}
\nu_k(n)=n\binom{k}{2}-3\tau_k(n).
\end{equation}
Petingi and Rodr\'iguez~\cite{PR} conjectured that, for every admissible $n$ and every integer $m\ge0$, $\nu_k(n+m(k+1))=\nu_k(n)$, and, more strongly, that every extremal graph on $n+m(k+1)$ vertices is the disjoint union of an extremal graph on $n$ vertices and $m$ copies of $K_{k+1}$. The conjecture is immediate for $k=1$, since every $1$-regular graph is a disjoint union of copies of $K_2$. Petingi and Rodr\'iguez~\cite{PR} verified the conjecture for $k=2,3$, and Gong, Yang, Wang and Chen~\cite{GYWC} settled the cases $k=4,5$.

The identity~\eqref{intro_path_minimum}, together with Corollary~\ref{cdk_consequences} and Theorems~\ref{maineven} and~\ref{mainperiodicity}, yields the following consequence for the conjecture of Petingi and Rodr\'iguez~\cite{PR}.
\begin{corollary}\label{consequences}
The conjecture of Petingi and Rodr\'iguez \cite{PR} holds whenever $k$ is odd, or $r$ is odd, or $3r\le2(k+1)$. Furthermore, for every even $k\ge18$, the conjecture fails in general, but holds for all $n\ge Ck$, where $C$ is an absolute constant.
\end{corollary}

\begin{remark}
For every even $k\ge18$ and every integer $m\ge0$, Corollary~\ref{cdk_consequences} rules out reducing an extremal graph of order $3k+m(k+1)$ to order $2k-1$ by deleting $m+1$ copies of $K_{k+1}$. The eventual recurrence in Corollary~\ref{consequences} applies to orders $n+m(k+1)$ with starting order $n\ge Ck$.
\end{remark}

From another perspective, triangle maximization in regular graphs is also equivalent to minimizing the number of independent sets of size three. Denote by $i_3(G)$ the number of sets of three vertices that are pairwise nonadjacent in $G$.
For a $k$-regular graph $G$ on $n$ vertices, Alexander and Mink~\cite{AM} give the identity
\begin{equation}\label{intro_independent_identity}
i_3(G)+\tau(G)=\binom{n}{3}-\frac{nk(n-1-k)}2.
\end{equation}
Independent set counting has been extensively studied~\cite{ACM,Alon,EG,Kahn,SSSZ,Zhao}. In particular, a result of Cutler and Radcliffe~\cite{CRtotal} implies that, when $k+1$ divides $n$, disjoint unions of copies of $K_{k+1}$ minimize the number of independent sets of size three among $k$-regular graphs on $n$ vertices. The relation $\tau(\ov G)=i_3(G)$ also connects triangle maximization in $k$-regular graphs to triangle minimization in regular graphs of degree $n-1-k$. Building on work of Lo~\cite{Lo}, He, Hou, Ma and Xie~\cite{HHMX} determined the minimum number of triangles in $d$-regular graphs on sufficiently large odd $n$ vertices when $d$ is even and $2n/5<d<n/2$. Using~\eqref{intro_independent_identity} and Theorems~\ref{mainodd} and~\ref{mainbounds}, we obtain the following lower bound on the number of independent sets of size three in $k$-regular graphs.
\begin{corollary}\label{independent_consequences}
Every $k$-regular graph $G$ on $n$ vertices satisfies
\[
i_3(G)\ge\binom n3-\frac{nk(n-1-k)}2-f(n,k),
\]
with equality if and only if $G$ has one of the forms in Theorem~\ref{mainodd}.
\end{corollary}


This paper is organized as follows. Section~\ref{sec:preliminaries} contains some preliminaries. In Section~\ref{sec:proof-mainodd}, we prove Theorem~\ref{mainodd} by adapting Chase's closed neighborhood counting argument~\cite{Chase}. In Section~\ref{sec:proof-mainbounds}, we prove Theorem~\ref{mainbounds} using an existence condition for triangle-free regular graphs and an explicit construction. In Section~\ref{sec:counterexamples}, we prove Theorem~\ref{maineven} using an explicit construction and a lower bound on triangle counts in complements of regular graphs. Finally, in Section~\ref{sec:eventual-periodicity}, we prove Theorem~\ref{mainperiodicity}.

\section{Preliminaries}\label{sec:preliminaries}
We use standard graph-theoretic notation~\cite{BM} and all graphs considered are finite and simple. For a positive integer $m$, write $[m]=\left\{1,\ldots,m\right\}$. Let $e(G)=|E(G)|$. For disjoint sets $S,T\subseteq V(G)$, let $e_G(S,T)$ denote the number of edges with one endpoint in $S$ and the other in $T$. For $v\in V(G)$, let $N_G(v)$ and $N_G[v]=N_G(v)\cup\left\{v\right\}$ denote its \emph{open neighborhood} and \emph{closed neighborhood}, respectively. The \emph{degree} of $v$ is denoted by $d_G(v)=|N_G(v)|$. For $S\subseteq V(G)$, let $G[S]$ be the \emph{induced subgraph} on $S$ and $G-S=G[V(G)\setminus S]$. A vertex $u$ is a \emph{nonneighbor} of a vertex $v$ if $u\ne v$ and $u$ is not adjacent to $v$. An unordered pair of distinct nonadjacent vertices is called a \emph{nonedge}. In the proofs we will repeatedly use the handshaking lemma, which gives $\sum_{v\in V(G)} d_G(v)=2e(G)$.

To simplify notation, we introduce the following two quantities. For a graph $G$ on $n$ vertices, let $$\Psi(G)=n+6e(G)+6\tau(G).$$
Note that the quantity $\Psi(G)$ can be interpreted as the number of ordered triples $(x,y,z)$ of vertices in $G$, with repetitions allowed, such that any two distinct vertices among $x,y,z$ are adjacent. Triples with one, two or three distinct vertices contribute $n$, $6e(G)$ or $6\tau(G)$, respectively.

For integers $k\ge1$ and $n\ge0$, write $n=q(k+1)+r$ with integers $q\ge0$ and $0\le r<k+1$, and define $$g_{k+1}(n)=q(k+1)^3+r^3.$$
Note that $\Psi$ is additive over disjoint unions and satisfies $\Psi(K_s)=s^3$ for every integer $s\ge0$. Thus $g_{k+1}(n)=\Psi(qK_{k+1}\cup K_r)$.

If $G_0$ is $k$-regular on $k+1+r$ vertices, its complement is $r$-regular. Using $i_3(G_0)+\tau(G_0)=\binom{k+1+r}{3}-kr(k+1+r)/2$ from~\eqref{intro_independent_identity}, together with $i_3(G_0)=\tau(\ov{G_0})$ and $2e(G_0)=k(k+1+r)$, we obtain from the definition of $\Psi$
\begin{equation}\label{nucomp}
\Psi(G_0)=(k+1)^3+r^3-6\tau(\ov{G_0}).
\end{equation}
For a $k$-regular graph $G$ on $n$ vertices, the handshaking lemma gives $e(G)=nk/2$. Hence
\begin{equation}\label{psireg}
\Psi(G)=n(3k+1)+6\tau(G).
\end{equation}

Caro and Tuza~\cite{CT} proved that, for every odd integer $n\ge1$, the maximum number of edges in a triangle-free regular graph on $n$ vertices is $n\lfloor n/5\rfloor$. Together with Petersen's $2$-factorization theorem~\cite{Petersen}, which states that every regular graph of positive even degree admits a decomposition into edge-disjoint spanning $2$-regular subgraphs, this yields the following lemma.
\begingroup
\begin{lemma}\label{trianglefreeexistence}
Let $n\ge1$ be odd and $d\ge0$ be even. A triangle-free $d$-regular graph on $n$ vertices exists if and only if $5d\le2n$.
\end{lemma}
\begin{proof}
The case $d=0$ is immediate, so assume $d>0$. By the result of Caro and Tuza~\cite{CT}, any triangle-free $d$-regular graph on $n$ vertices satisfies $nd/2\le n\lfloor n/5\rfloor$, and hence $5d\le2n$. Conversely, assume that $5d\le2n$. Since $d$ is even, we have $d/2\le\lfloor n/5\rfloor$. The result of Caro and Tuza~\cite{CT} also guarantees that the maximum edge count $n\lfloor n/5\rfloor$ is attained by a triangle-free regular graph on $n$ vertices, which is necessarily $2\lfloor n/5\rfloor$-regular by the handshaking lemma. By Petersen's theorem~\cite{Petersen}, we can decompose this graph into $2$-factors and take the union of any $d/2$ of them to obtain the required triangle-free $d$-regular graph.
\end{proof}
\endgroup

We will also use the following construction.
\begin{construction}\label{construction:regular-bipartite}
Let $s,t$ be integers with $1\le t\le s$. Take two disjoint vertex sets $\left\{x_0,\ldots,x_{s-1}\right\}$ and $\left\{y_0,\ldots,y_{s-1}\right\}$, and join $x_i$ to $y_{i+j}$ for $0\le i<s$ and $0\le j<t$, with subscripts taken modulo $s$. The resulting graph is simple, $t$-regular and bipartite, and the edges $x_i y_i$ for $0\le i<s$ form a perfect matching.
\end{construction}

\section{Proof of Theorem~\ref{mainodd}}\label{sec:proof-mainodd}
For a graph $G$, let $\Delta(G)$ denote its maximum degree and $c_{1,3}(G)$ the number of induced copies of $K_{1,3}$ in $G$. In Lemma~\ref{closedneighborhood}, we adapt Chase's quadruple counting argument~\cite{Chase}, using closed neighborhoods and allowing repeated vertices, to establish a counting identity for $\Psi$ that will be used to derive a sharp upper bound under a maximum degree constraint and characterize the equality cases.

\begin{lemma}\label{closedneighborhood}
Let $G$ be a graph. Then
\[
\sum_{v\in V(G)}\bigl(\Psi(G)-\Psi(G-N_G[v])\bigr)
=\sum_{v\in V(G)}(d_G(v)+1)^3-6c_{1,3}(G).
\]
\end{lemma}
\begin{proof}
Let $\mathcal S$ consist of ordered quadruples $(v,x,y,z)\in V(G)^4$ with $x,y,z\in N_G[v]$. Let $\mathcal T$ consist of ordered quadruples $(v,x,y,z)\in V(G)^4$ such that any two distinct vertices among $x,y,z$ are adjacent and at least one of $x,y,z$ belongs to $N_G[v]$. Repeated vertices are allowed in both $\mathcal S$ and $\mathcal T$. Using the counting interpretation of $\Psi$, we obtain
\[
|\mathcal S|=\sum_{v\in V(G)}(d_G(v)+1)^3,\qquad
|\mathcal T|=\sum_{v\in V(G)}\bigl(\Psi(G)-\Psi(G-N_G[v])\bigr).
\]
Fix $U\subseteq V(G)$ with $1\le |U|\le4$, and count quadruples satisfying $\left\{v,x,y,z\right\}=U$. The nonzero counts are listed in Table~\ref{counts}.

\begin{table}[htbp]
\centering
\caption{Counts of ordered quadruples on a fixed vertex set $U$.}\label{counts}
\begin{tabular}{c*{8}{c}}
\toprule
$G[U]$&$K_1$&$K_2$&$P_3$&$K_3$&$K_{1,3}$&paw&diamond&$K_4$\\
\midrule
$\mathcal S$&1&14&12&36&6&6&12&24\\
$\mathcal T$&1&14&12&36&0&6&12&24\\
\bottomrule
\end{tabular}
\par\smallskip
\begin{minipage}{0.94\linewidth}\small
\emph{Note.} A \emph{paw} is a triangle with a pendant vertex, and a \emph{diamond} is $K_4$ with one edge deleted. All other isomorphism types contribute zero to both counts.
\end{minipage}
\end{table}

For $K_1$, $K_2$ and $K_3$, every ordered quadruple using all vertices belongs to both families, giving $1$, $2^4-2=14$ and $3^4-3\cdot2^4+3=36$, respectively. For $P_3$, a quadruple in $\mathcal S$ must have $v$ as the middle vertex and both endpoints among $x,y,z$, giving $3^3-2\cdot2^3+1=12$ choices. In $\mathcal T$, $v$ must be an endpoint, and $(x,y,z)$ must use each of the other two vertices at least once, giving $2(2^3-2)=12$ choices.

When $|U|=4$, all vertices are distinct. In $\mathcal S$, $v$ is adjacent to the other three vertices. In $\mathcal T$, $x,y,z$ form a triangle and at least one is adjacent to $v$. For $K_{1,3}$, the paw, the diamond and $K_4$, there are respectively $1,1,2,4$ choices of $v$ in $\mathcal S$ and $0,1,2,4$ in $\mathcal T$, followed by $3!=6$ orders of $x,y,z$. Clearly, disconnected graphs, $P_4$ and the cycle on four vertices contribute zero to both families. Hence $U$ contributes $6$ to $|\mathcal S|-|\mathcal T|$ exactly when it induces a copy of $K_{1,3}$, and zero otherwise. This proves Lemma \ref{closedneighborhood}.
\end{proof}

We now use Lemma~\ref{closedneighborhood} to bound $\Psi(G)$ and characterize the cases of equality.

\begin{lemma}\label{psilemma}
Let $n=q(k+1)+r$ for integers $k,q,r$ with $k\ge1$, $q\ge0$ and $0\le r<k+1$. Every graph $G$ on $n$ vertices with $\Delta(G)\le k$ satisfies
\[ \Psi(G)\le g_{k+1}(n). \]
Equality holds if and only if one of the following holds:
\begin{enumerate}[label=(\roman*),leftmargin=2em]
\item [$(i).$] $G\cong qK_{k+1}\cup K_r$.
\item [$(ii).$] $G\cong(q-1)K_{k+1}\cup G_0$, where $q,r\ge1$ and $G_0$ is a $k$-regular graph on $k+1+r$ vertices whose complement is triangle-free.
\end{enumerate}
\end{lemma}

\begin{proof}
For fixed $k$, we first prove the inequality $\Psi(G)\le g_{k+1}(n)$ by induction on $n$. The null graph gives the base case. For $n>0$, Lemma~\ref{closedneighborhood} gives a vertex $v\in V(G)$ such that
\begin{equation}\label{vertexbound}
\Psi(G)-\Psi(G-N_G[v])\le(d_G(v)+1)^3.
\end{equation}
For every vertex $v$ of $G$, since $1\le d_G(v)+1\le\min\{k+1,n\}$, the definition of $g_{k+1}$ gives the following identity
\begin{equation*}
\begin{aligned}
g_{k+1}(n)-g_{k+1}\bigl(n-d_G(v)-1\bigr)-(d_G(v)+1)^3
=\begin{cases}
3r(d_G(v)+1)(r-d_G(v)-1),&d_G(v)<r,\\
3(k+1+r)(d_G(v)+1-r)(k-d_G(v)),&r\le d_G(v)\le k.
\end{cases}
\end{aligned}
\end{equation*}
Note that the right-hand side of this identity is nonnegative, with equality if and only if $d_G(v)=r-1$ or $d_G(v)=k$. Since $G-N_G[v]$ has $n-d_G(v)-1<n$ vertices and maximum degree at most $k$, the induction hypothesis, together with~\eqref{vertexbound} and the preceding identity, yields
\begin{align*}
\Psi(G)
\le\Psi(G-N_G[v])+(d_G(v)+1)^3
\le g_{k+1}\bigl(n-d_G(v)-1\bigr)+(d_G(v)+1)^3
\le g_{k+1}(n).
\end{align*}
This proves the inequality in Lemma~\ref{psilemma}.

We now prove that equality holds if and only if $G$ has one of the forms stated in Lemma~\ref{psilemma}. Assume that $\Psi(G)=g_{k+1}(n)$. The null graph has the form described in Lemma~\ref{psilemma}(i), so we may assume $n>0$. Applying the inequality in Lemma~\ref{psilemma} to $G-N_G[v]$ and using the nonnegativity of the difference computed above, we obtain, for every vertex $v\in V(G)$,
\begin{equation}\label{equalitychain}
\Psi(G)-\Psi(G-N_G[v])
\ge g_{k+1}(n)-g_{k+1}\bigl(n-(d_G(v)+1)\bigr)
\ge(d_G(v)+1)^3.
\end{equation}
By Lemma~\ref{closedneighborhood}, it holds that
\[ \sum_{v\in V(G)}\bigl(\Psi(G)-\Psi(G-N_G[v])-(d_G(v)+1)^3\bigr)=-6c_{1,3}(G). \]
By~\eqref{equalitychain}, every summand on the left is nonnegative, whereas the right-hand side is nonpositive. Thus $c_{1,3}(G)=0$ and $\Psi(G)-\Psi(G-N_G[v])-(d_G(v)+1)^3=0$. Hence equality holds throughout the inequality chain~\eqref{equalitychain} for every $v\in V(G)$, yielding
\begin{gather}
\Psi(G-N_G[v])=g_{k+1}\bigl(n-(d_G(v)+1)\bigr), \text{and}\label{hereditary}\\[6pt]
g_{k+1}(n)-g_{k+1}\bigl(n-(d_G(v)+1)\bigr)=(d_G(v)+1)^3.\label{degreeequality}
\end{gather}
Equation~\eqref{degreeequality} and the equality condition for the nonnegative difference above show that every vertex of $G$ has degree $k$ if $r=0$, and degree either $r-1$ or $k$ if $r>0$.

\medskip\noindent\textbf{Case 1.} $r=0$.\par 
In this case, $G$ is $k$-regular and $\Psi(G)=g_{k+1}(n)=n(k+1)^2$, so the identity $\Psi(G)=n(3k+1)+6\tau(G)$ in~\eqref{psireg} gives $3\tau(G)=nk(k-1)/2$. Clearly, each vertex of $G$ belongs to at most $k(k-1)/2$ triangles. Since the sum of these vertex counts is $3\tau(G)=nk(k-1)/2$, every vertex attains this bound. Hence $G[N_G[v]]\cong K_{k+1}$ for every $v\in V(G)$. Since $G$ is $k$-regular, each such complete graph is a component of $G$, and therefore $G\cong qK_{k+1}$.

\medskip\noindent\textbf{Case 2.} $r>0$ and $G$ has a vertex of degree less than $k$.\par 
We first choose $v\in V(G)$ with $d_G(v)=r-1$. Recall that $n=q(k+1)+r$. By~\eqref{hereditary}, $G-N_G[v]$ attains equality at order $n-r=q(k+1)$.  For $q\ge1$, Case~1 gives $G-N_G[v]\cong qK_{k+1}$, which also holds trivially when $q=0$. The condition $\Delta(G)\le k$ implies that no edge joins these complete components to $N_G[v]$. The $r$ vertices of $N_G[v]$ therefore have degree at most $r-1<k$, so each has degree exactly $r-1$. Thus $G[N_G[v]]\cong K_r$, giving the form described in Lemma~\ref{psilemma}(i).

\medskip\noindent\textbf{Case 3.} $r>0$ and $G$ is $k$-regular.\par 
In this case, we have $q\ge1$. Since $r>0$ and $G$ is $k$-regular, $G$ is not a disjoint union of complete graphs. Hence $G$ contains an induced path $vwu$ with $vu\notin E(G)$. In $G-N_G[v]$, the vertex $u$ has lost its neighbor $w$ and therefore has degree less than $k$. Equation~\eqref{hereditary} shows that $G-N_G[v]$ attains equality at order $n-(k+1)=(q-1)(k+1)+r$. Applying Case~2 with $q$ replaced by $q-1$ gives $G-N_G[v]\cong(q-1)K_{k+1}\cup K_r$. The $q-1$ copies of $K_{k+1}$ are also components of $G$. Deleting them leaves a $k$-regular graph $G_0$ on $k+1+r$ vertices. By additivity of $\Psi$ over disjoint unions, $\Psi(G_0)=\Psi(G)-(q-1)(k+1)^3=(k+1)^3+r^3$. The identity $\Psi(G_0)=(k+1)^3+r^3-6\tau(\ov{G_0})$ in~\eqref{nucomp} now gives $\tau(\ov{G_0})=0$, yielding the form described in Lemma~\ref{psilemma}(ii).

Conversely, the form described in Lemma~\ref{psilemma}(i) gives equality since $\Psi(K_t)=t^3$ and $\Psi$ is additive. For the form described in Lemma~\ref{psilemma}(ii), equation~\eqref{nucomp} gives $\Psi(G_0)=(k+1)^3+r^3$, so additivity again yields equality. This completes the proof.
\end{proof}

\begin{proof}[\textbf{Proof of Theorem~\ref{mainodd}}]
Write $n=q(k+1)+r$ as in Theorem~\ref{mainodd}. For every $k$-regular graph $G$ on $n$ vertices, the identity $\Psi(G)=n(3k+1)+6\tau(G)$ in~\eqref{psireg} and Lemma~\ref{psilemma} give
\begin{align}
\tau(G)=\frac{\Psi(G)-n(3k+1)}6
\le\frac{g_{k+1}(n)-n(3k+1)}6=f(n,k).\label{oddvalue}
\end{align}
If $r=0$, the graph $qK_{k+1}$ attains this bound. Henceforth assume $r>0$. 
	
If either $k$ or $r$ is odd, then $k+1+r$ is even since $nk$ is even. Let $h=(k+1+r)/2$. Since $1\le r<h$, taking $s:=h$ and $t:=r$ in Construction~\ref{construction:regular-bipartite} gives an $r$-regular bipartite graph $H$ on $k+1+r$ vertices, which is triangle-free. Note that $\ov H$ is $k$-regular on $k+1+r$ vertices. By the equality statement of Lemma~\ref{psilemma}, the graph $(q-1)K_{k+1}\cup\ov H$ attains equality in~\eqref{oddvalue}. Now assume that $k$ and $r$ are even and $3r\le2(k+1)$. Then $k+1+r$ is odd and $5r\le2(k+1+r)$. By Lemma~\ref{trianglefreeexistence}, there is a triangle-free $r$-regular graph $H$ on $k+1+r$ vertices. Again, the graph $(q-1)K_{k+1}\cup\ov H$ attains equality in~\eqref{oddvalue}. Thus $\tau_k(n)=f(n,k)$, proving the formula in Theorem~\ref{mainodd}.

According to equation~\eqref{oddvalue}, a $k$-regular graph $G$ on $n$ vertices is extremal  if and only if it attains equality in Lemma~\ref{psilemma}. If $r>0$, the first form in that lemma is not $k$-regular, and the second form is exactly $(q-1)K_{k+1}\cup\ov H$, where $H$ is any triangle-free $r$-regular graph on $k+1+r$ vertices. If $r=0$, the equality characterization gives $qK_{k+1}$. This completes the proof.
\end{proof}

\section{Proof of Theorem~\ref{mainbounds}}\label{sec:proof-mainbounds}
\begingroup
\begin{proof}[\textbf{Proof of Theorem~\ref{mainbounds}}]
Assume that $k$ and $r$ are even and $3r>2(k+1)$. We first prove the upper bound. Let $G$ be any $k$-regular graph on $n$ vertices. The identity $\Psi(G)=n(3k+1)+6\tau(G)$ in~\eqref{psireg}, together with Lemma~\ref{psilemma}, gives $\tau(G)\le f(n,k)$. Suppose that $\tau(G)=f(n,k)$. Since $r>0$, the graph $qK_{k+1}\cup K_r$ is not $k$-regular. Hence the equality statement of Lemma~\ref{psilemma} gives $G\cong(q-1)K_{k+1}\cup\ov H$, where $H$ is a triangle-free $r$-regular graph on $k+1+r$ vertices. Since $k+1+r$ is odd, Lemma~\ref{trianglefreeexistence} gives $5r\le2(k+1+r)$, and hence $3r\le2(k+1)$, a contradiction. Thus $\tau(G)<f(n,k)$. Moreover, $f(n,k)=q\binom{k+1}{3}+\binom r3-r(k+1-r)/2$ is an integer since $r$ is even. Hence $\tau_k(n)\le f(n,k)-1$, proving the upper bound.

For the lower bound, we construct an $r$-regular graph $H$ on $k+1+r$ vertices with exactly $r(2r-k-2)/4$ triangles (see Figure~\ref{fig:mainbounds-construction}). Let $A_1,B_1,A_2,B_2$ be disjoint independent sets with $|A_1|=|B_1|=r/2$ and $|A_2|=|B_2|=k/2$. Join every vertex of $A_1$ to every vertex of $B_2$, and every vertex of $A_2$ to every vertex of $B_1$. Since $3r>2(k+1)$ and $r\le k$ imply $1\le r-k/2-1\le r/2\le k/2$, taking $s:=r/2$ and $t:=r-k/2-1$ in Construction~\ref{construction:regular-bipartite} gives an $(r-k/2-1)$-regular bipartite graph between $A_1$ and $B_1$. Taking $s:=k/2$ and $t:=r/2$ in the same construction gives an $(r/2)$-regular bipartite graph between $A_2$ and $B_2$. Finally, add a vertex $z$ adjacent to every vertex of $A_1\cup B_1$. Denote the resulting graph by $H$.

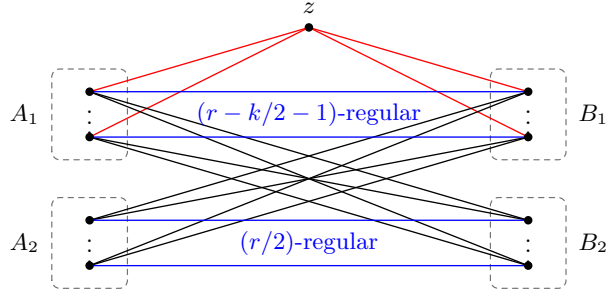
\begin{figure}[htbp]
\centering
\begin{tikzpicture}[
    every node/.style={font=\small},
    line width=.5pt,line cap=round,line join=round,
    groupbox/.style={draw=black!55,dash pattern=on 2pt off 2pt,
        rounded corners=3pt,line width=.45pt}
]
\draw[groupbox] (-3.4,.25) rectangle (-2.4,1.45);
\draw[groupbox] (2.4,.25) rectangle (3.4,1.45);
\draw[groupbox] (-3.4,-1.45) rectangle (-2.4,-.25);
\draw[groupbox] (2.4,-1.45) rectangle (3.4,-.25);
\coordinate (A1a) at (-2.9,1.15);
\coordinate (A1b) at (-2.9,.55);
\coordinate (B1a) at (2.9,1.15);
\coordinate (B1b) at (2.9,.55);
\coordinate (A2a) at (-2.9,-1.15);
\coordinate (A2b) at (-2.9,-.55);
\coordinate (B2a) at (2.9,-1.15);
\coordinate (B2b) at (2.9,-.55);
\coordinate (z) at (0,2.0);
\foreach \v in {A1a,A1b,B1a,B1b}{\draw[red] (z)--(\v);}
\foreach \u in {A1a,A1b}{\foreach \v in {B2a,B2b}{\draw (\u)--(\v);}}
\foreach \u in {B1a,B1b}{\foreach \v in {A2a,A2b}{\draw (\u)--(\v);}}
\draw[blue] (A1a)--(B1a);
\draw[blue] (A1b)--(B1b);
\draw[blue] (B2a)--(A2a);
\draw[blue] (B2b)--(A2b);
\foreach \v in {z,A1a,A1b,B1a,B1b,B2a,B2b,A2a,A2b}{\node[vertex] at (\v) {};}
\node[above=2pt] at (z) {$z$};
\node at (-2.9,.85) {$\vdots$};
\node at (2.9,.85) {$\vdots$};
\node at (-2.9,-.85) {$\vdots$};
\node at (2.9,-.85) {$\vdots$};
\node[left=4pt] at (-3.3,.85) {$A_1$};
\node[right=4pt] at (3.3,.85) {$B_1$};
\node[left=4pt] at (-3.3,-.85) {$A_2$};
\node[right=4pt] at (3.3,-.85) {$B_2$};
\node[inner sep=1pt,text=blue] at (0,.85) {$(r-k/2-1)$-regular};
\node[inner sep=1pt,text=blue] at (0,-.85) {$(r/2)$-regular};
\end{tikzpicture}
\caption{An illustration of $H$.}\label{fig:mainbounds-construction}
\end{figure}

Note that every vertex of $A_1\cup B_1$ has degree $k/2+(r-k/2-1)+1=r$, every vertex of $A_2\cup B_2$ has degree $r/2+r/2=r$, and $z$ has degree $r$. Thus $H$ is $r$-regular on $k+1+r$ vertices. Since $H-z$ is bipartite with parts $A_1\cup A_2$ and $B_1\cup B_2$, every triangle of $H$ contains $z$. Therefore, $\tau(H)=(r/2)(r-k/2-1)=r(2r-k-2)/4$. Let $G=(q-1)K_{k+1}\cup\ov H$. Since $\ov H$ is $k$-regular on $k+1+r$ vertices, $G$ is $k$-regular on $n$ vertices. Using $\Psi(\ov H)=(k+1)^3+r^3-6\tau(H)$ from~\eqref{nucomp} and the additivity of $\Psi$ over disjoint unions, we obtain $\Psi(G)=(q-1)(k+1)^3+\Psi(\ov H)=q(k+1)^3+r^3-6\tau(H)$. The identity $\Psi(G)=n(3k+1)+6\tau(G)$ in~\eqref{psireg} then gives
\[
\tau(G)=f(n,k)-\tau(H)=f(n,k)-\frac{r(2r-k-2)}4.
\]
Since $\tau_k(n)\ge\tau(G)$, the lower bound follows.
\end{proof}
\endgroup

\section{Proof of Theorem~\ref{maineven}}\label{sec:counterexamples}
We prove the two bounds in Theorem~\ref{maineven} separately. Construction~\ref{construction} gives the lower bound at order $3k$, while Lemma~\ref{trianglelemma}, applied to the complement, gives the upper bound at order $2k-1$.

\begin{construction}[$G_k$]\label{construction}
Let $k\ge2$ be even. Take three disjoint copies of $K_k$ with vertex sets $V_1,V_2,V_3$. Partition each $V_i$ into two sets $V_{ij}$, indexed by $j\in\left\{1,2,3\right\}\setminus\left\{i\right\}$, each of size $k/2$. For each $1\le i<j\le3$, add a perfect matching between $V_{ij}$ and $V_{ji}$.
\end{construction}

Each vertex $z\in V_i$ has $k-1$ neighbors in $V_i$ and exactly one neighbor outside $V_i$, so $G_k$ is $k$-regular. Note  that the added edges form a matching on $V(G_k)$, and no triangle can contain any of them. Thus every triangle lies in one of the three copies of $K_k$, giving $\tau(G_k)=3\binom{k}{3}$. 

\begin{lemma}\label{trianglelemma}
Let $d\ge16$ be even. Every $d$-regular graph $G$ on $2d+3$ vertices satisfies $\tau(G)\ge3d/2+3$.
\end{lemma}
\begin{proof}
Suppose that there exists such a graph $G$ with $\tau(G)\le3d/2+2$. Note that each triangle contributes three to $\sum_{z\in V(G)} e(G[N_G(z)])$. By the supposition, the average of $e(G[N_G(z)])$ over all vertices $z\in V(G)$ satisfies
\[
\frac{1}{2d+3}\sum_{z\in V(G)}e(G[N_G(z)])
=\frac{3\tau(G)}{2d+3}
\le\frac{3(3d/2+2)}{2d+3}<3.
\]
Since these edge counts are integers, there is a vertex $u\in V(G)$ with $\ell:=e(G[N_G(u)])\le2$. Let $A=N_G(u)$, $B=V(G)\setminus N_G[u]$ and $H=G[B]$; see Figure~\ref{figpartition}.
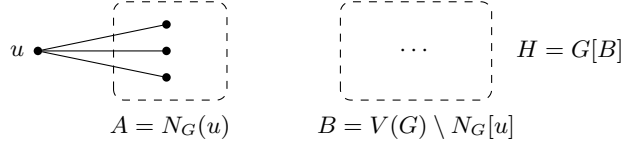
\begin{figure}[htbp]
\centering
\begin{tikzpicture}[every node/.style={font=\small}]
\node[vertex,label=left:$u$] (u) at (0,0) {};
\draw[dashed,rounded corners] (1,-.65) rectangle (2.5,.65);
\draw[dashed,rounded corners] (4,-.65) rectangle (6,.65);
\foreach \y in {-.35,0,.35}{\draw (u)--(1.7,\y);\node[vertex] at (1.7,\y) {};}
\node at (5,0){$\cdots$};\node at (1.75,-1){$A=N_G(u)$};
\node at (5,-1){$B=V(G)\setminus N_G[u]$};\node[right] at (6.2,0){$H=G[B]$};
\end{tikzpicture}
\caption{An illustration of the partition $V(G)=\left\{u\right\}\cup A\cup B$.}\label{figpartition}
\end{figure}

Note that  $|A|=d$, $|B|=d+2$ and $e(G[A])=\ell$. Every $v\in A$ has $d-1-d_{G[A]}(v)$ neighbors and hence $3+d_{G[A]}(v)$ nonneighbors in $B$. Every $w\in B$ has $d-d_H(w)$ neighbors and hence $d_H(w)$ nonneighbors in $A$. Counting nonedges with one endpoint in each of $A$ and $B$ gives
\[
\sum_{w\in B}d_H(w)=\sum_{v\in A}(3+d_{G[A]}(v))=3d+2\ell.
\]
Thus the handshaking lemma yields
\begin{equation}\label{edgesH}
e(H)=\frac{3d}{2}+\ell.
\end{equation}

For a real number $\alpha$, let $\pp{\alpha}=\max\left\{\alpha,0\right\}$, and let
\[
\beta=\sum_{ww'\in E(H)}\pp{d-d_H(w)-d_H(w')}.
\]
For every edge $ww'\in E(H)$, its endpoints have at least $\pp{d-d_H(w)-d_H(w')}$ common neighbors in $A$. Thus $\beta$ is a lower bound for the number of triangles with one vertex in $A$ and two in $B$.

Since $\ell\le2$, there is no triangle entirely in $A$. Hence each triangle of $G$ contains $u$, meets both $A$ and $B$, or lies entirely in $B$. Therefore,
\begin{align*}
\tau(G)={}&\ell
+\sum_{ww'\in E(H)}|N_G(w)\cap N_G(w')\cap A|
+\sum_{vv'\in E(G[A])}|N_G(v)\cap N_G(v')\cap B|+\tau(H)\\[6pt]
\ge{}&\ell+\beta
+\sum_{vv'\in E(G[A])}\bigl[(d+2)-(3+d_{G[A]}(v))-(3+d_{G[A]}(v'))\bigr]+\tau(H)\\[6pt]
={}&\ell+\beta
+\sum_{vv'\in E(G[A])}\bigl(d-4-d_{G[A]}(v)-d_{G[A]}(v')\bigr)+\tau(H).
\end{align*}
Considering the possible values of $\ell$ and the corresponding structures of $G[A]$, we obtain
\begin{equation}\label{trianglecases}
\raisebox{\dimexpr\depth-1.2\dp\strutbox\relax}{$\displaystyle
\tau(G)\ge\tau(H)+\beta+
\begin{cases}
0,&\ell=0,\\
d-5,&\ell=1,\\
2d-10,&\ell=2\text{ and }\Delta(G[A])=1,\\
2d-12,&\ell=2\text{ and }\Delta(G[A])=2.
\end{cases}
$}
\end{equation}

We next derive a lower bound for $\tau(H)+\beta$. Let $w_0$ be a vertex of maximum degree in $H$. We claim that $d/2\le d_H(w_0)\le d$. Clearly, $d_H(w_0)\le d_G(w_0)\leq d$. If $d_H(w_0)\le d/2-1$, every summand in the definition of $\beta$ is at least $2$, so~\eqref{trianglecases} and~\eqref{edgesH} imply $\tau(G)\ge\beta\ge2e(H)=3d+2\ell>3d/2+2$, contrary to the supposition. This proves the claim.

Let $H_0=H-\left\{w_0\right\}$. Since $N_H(w_0)\subseteq V(H_0)$, the handshaking lemma gives
\begin{equation}\label{Sdegreebound}
\sum_{z\in N_H(w_0)}d_{H_0}(z)\le\sum_{z\in V(H_0)}d_{H_0}(z)=2e(H_0).
\end{equation}
Using $e(H)=3d/2+\ell$ from~\eqref{edgesH}, we obtain $e(H_0)=e(H)-d_H(w_0)=3d/2+\ell-d_H(w_0)$. Note that every edge of $H_0[N_H(w_0)]$ forms a triangle with $w_0$. Also, $\sum_{z\in N_H(w_0)}d_{H_0}(z)$ counts each edge of $H_0$ with both endpoints in $N_H(w_0)$ twice and each edge with exactly one endpoint there once. Thus $\sum_{z\in N_H(w_0)}d_{H_0}(z)-e(H_0[N_H(w_0)])\le e(H_0)$. Consequently,
\begin{equation}\label{Sbound}
\sum_{z\in N_H(w_0)}d_{H_0}(z)\le e(H_0)+e(H_0[N_H(w_0)])\le e(H_0)+\tau(H).
\end{equation}
Using the handshaking lemma and the fact that any two distinct edges of a simple graph have at most one common endpoint, we obtain
\begin{equation}\label{degreesquares}
\begin{aligned}
\sum_{z\in V(H_0)}d_{H_0}(z)^2
=2\sum_{z\in V(H_0)}\binom{d_{H_0}(z)}2+2e(H_0)
\le2\binom{e(H_0)}2+2e(H_0)=e(H_0)^2+e(H_0).
\end{aligned}
\end{equation}
For $z\in N_H(w_0)$, we have $d_H(z)=d_{H_0}(z)+1$. For $z\in V(H_0)\setminus N_H(w_0)$, we have $d_H(z)=d_{H_0}(z)$. Since each vertex $z\in V(H_0)$ is an endpoint of exactly $d_{H_0}(z)$ edges of $H_0$, we obtain
\begin{equation}\label{edgesum}
\begin{aligned}
\sum_{ww'\in E(H_0)}(d_H(w)+d_H(w'))
=\sum_{z\in V(H_0)}d_H(z)d_{H_0}(z)
=\sum_{z\in V(H_0)}d_{H_0}(z)^2+\sum_{z\in N_H(w_0)}d_{H_0}(z).
\end{aligned}
\end{equation}

Split the sum for $\beta$ according to whether an edge is incident with $w_0$ or lies in $H_0$. Since $d_H(z)=d_{H_0}(z)+1$ for $z\in N_H(w_0)$, we have
\[
\sum_{z\in N_H(w_0)}d_H(z)
=\sum_{z\in N_H(w_0)}\bigl(d_{H_0}(z)+1\bigr)
=\sum_{z\in N_H(w_0)}d_{H_0}(z)+d_H(w_0).
\]
Using $\sum_{i=1}^{m}\pp{a_i}\ge\pp{\sum_{i=1}^{m}a_i}$ for real $a_1,\ldots,a_m$, together with~\eqref{Sdegreebound}, \eqref{degreesquares} and~\eqref{edgesum}, we obtain
\begin{align}
\beta
&=\sum_{z\in N_H(w_0)}\pp{d-d_H(w_0)-d_H(z)}
+\sum_{ww'\in E(H_0)}\pp{d-d_H(w)-d_H(w')}\notag\\[6pt]
&\ge\pp{d_H(w_0)(d-d_H(w_0))-\sum_{z\in N_H(w_0)}d_H(z)}
+\pp{de(H_0)-\sum_{ww'\in E(H_0)}(d_H(w)+d_H(w'))}\notag\\[8pt]
&=\pp{d_H(w_0)(d-d_H(w_0)-1)-\sum_{z\in N_H(w_0)}d_{H_0}(z)}
+\pp{de(H_0)-\sum_{z\in V(H_0)}d_{H_0}(z)^2-\sum_{z\in N_H(w_0)}d_{H_0}(z)}\notag\\[6pt]
&\ge\pp{d_H(w_0)(d-d_H(w_0)-1)-\sum_{z\in N_H(w_0)}d_{H_0}(z)}+\pp{de(H_0)-e(H_0)^2-3e(H_0)}.\label{betabound}
\end{align}
We combine $\tau(H)\ge\pp{\sum_{z\in N_H(w_0)}d_{H_0}(z)-e(H_0)}$, which follows from~\eqref{Sbound}, with the lower bound in~\eqref{betabound}. Using $\pp{a}+\pp{b}\ge\pp{a+b}$ for real $a,b$ and $e(H_0)=3d/2+\ell-d_H(w_0)$, as obtained from~\eqref{edgesH}, we then obtain
\begingroup
\medmuskip=3mu
\begin{align}\notag
\tau(H)+\beta
&\ge\pp{\sum_{z\in N_H(w_0)}d_{H_0}(z)-e(H_0)}+\pp{d_H(w_0)(d-d_H(w_0)-1)-\sum_{z\in N_H(w_0)}d_{H_0}(z)}\notag\\[6pt]
&\quad +\pp{de(H_0)-e(H_0)^2-3e(H_0)}\notag\\[6pt]
&\ge\pp{d_H(w_0)(d-d_H(w_0)-1)-e(H_0)}+\pp{de(H_0)-e(H_0)^2-3e(H_0)}\notag\\[6pt]
&=\pp{d_H(w_0)(d-d_H(w_0))-\frac{3d}{2}-\ell}
+\pp{\left(\frac{3d}{2}+\ell-d_H(w_0)\right)
\left(d_H(w_0)-\frac d2-\ell-3\right)}.\label{keybound}
\end{align}
\endgroup

\begin{claim}\label{clm:triangle-complement-bound}
$\tau(H)+\beta\ge d^2/4-3d/2-\ell^2-7\ell-9$.
\end{claim}
\begin{proof}
Suppose for a contradiction that $\tau(H)+\beta<d^2/4-3d/2-\ell^2-7\ell-9$. By~\eqref{keybound}, we have
\begin{align*}
\tau(H)+\beta
&\ge\pp{d_H(w_0)(d-d_H(w_0))-\frac{3d}{2}-\ell}
\ge d_H(w_0)(d-d_H(w_0))-\frac{3d}{2}-\ell.
\end{align*}
Combining this with the supposition, we obtain
\begin{equation}\label{firstproduct}
\left(\frac d2+\ell+3-d_H(w_0)\right)\left(d_H(w_0)-\frac d2+\ell+3\right)<0.
\end{equation}
Since $d_H(w_0)\ge d/2$, $d_H(w_0)-d/2+\ell+3>0$, and hence, by~\eqref{firstproduct}, $d_H(w_0)>d/2+\ell+3$. Estimating the second positive-part term in~\eqref{keybound} in the same way and combining both estimates with the supposition, we obtain
\begin{align*}
0>\tau(H)+\beta-\left(\frac{d^2}{4}-\frac{3d}{2}-\ell^2-7\ell-9\right)\ge\left(d_H(w_0)-\frac d2-\ell-3\right)(2d-2d_H(w_0)-3).
\end{align*}
Since $d_H(w_0)-d/2-\ell-3>0$, it holds that $d_H(w_0)>d-3/2$. Since $d_H(w_0)$ is an integer with $d_H(w_0)\le d$, it follows that $d_H(w_0)\in\left\{d-1,d\right\}$. For these two values, the product $\left(3d/2+\ell-d_H(w_0)\right) \left(d_H(w_0)-d/2-\ell-3\right)$ in~\eqref{keybound} takes the respective values
\[
\begin{aligned}
\left(\frac d2+\ell+1\right)\left(\frac d2-\ell-4\right)
=\frac{d^2}{4}-\frac{3d}{2}-\ell^2-5\ell-4 \quad\text{and}\quad
\left(\frac d2+\ell\right)\left(\frac d2-\ell-3\right)
=\frac{d^2}{4}-\frac{3d}{2}-\ell^2-3\ell.
\end{aligned}
\]
These values exceed the claimed lower bound by $2\ell+5$ and $4\ell+9$, respectively, contradicting the supposition and proving the claim.
\end{proof}
\medskip
Combining Claim~\ref{clm:triangle-complement-bound} with~\eqref{trianglecases}, using $2d-12$ for both configurations with $\ell=2$, gives
\[
\tau(G)-\left(\frac{3d}{2}+3\right)\ge
\begin{cases}
\frac{d^2}{4}-3d-12,&\ell=0,\\
\frac{d^2}{4}-2d-25,&\ell=1,\\
\frac{d^2}{4}-d-42,&\ell=2.
\end{cases}
\]
At $d=16$, these expressions are $4,7,6$, respectively, and each is increasing for $d\ge16$. This contradicts the initial supposition and proves Lemma~\ref{trianglelemma}.
\end{proof}

\begin{proof}[\textbf{Proof of Theorem~\ref{maineven}}]
Construction~\ref{construction} gives a $k$-regular graph $G_k$ on $3k$ vertices with $\tau(G_k)=3\binom{k}{3}=f(3k,k)-3k/2+1$, proving the first bound. Let $G_0$ be any $k$-regular graph on $2k-1$ vertices. Then $\ov{G_0}$ is $(k-2)$-regular on $2(k-2)+3$ vertices. Applying Lemma~\ref{trianglelemma} to $\ov{G_0}$ with $d=k-2$, we obtain $\tau(\ov{G_0})\ge3k/2$. Using $i_3(G_0)+\tau(G_0)=\binom{2k-1}{3}-(2k-1)k(k-2)/2$ from~\eqref{intro_independent_identity}, together with $i_3(G_0)=\tau(\ov{G_0})$, we obtain
\begin{align*}
\tau(G_0)=\binom{2k-1}{3}-\frac{(2k-1)k(k-2)}2-\tau(\ov{G_0})
\le\binom{2k-1}{3}-\frac{(2k-1)k(k-2)}2-\frac{3k}{2}
=f(2k-1,k)-\frac{3k}{2}.
\end{align*}
This proves the second bound. Since $f(3k,k)-f(2k-1,k)=\binom{k+1}{3}$, the two bounds give $\tau_k(3k)\ge\tau_k(2k-1)+\binom{k+1}{3}+1$, proving the strict inequality in the theorem.
\end{proof}

\begingroup
\section{Proof of Theorem~\ref{mainperiodicity}}\label{sec:eventual-periodicity}
Throughout this section, let $G$ be an extremal $k$-regular graph on $n=q(k+1)+r$ vertices, where $k\ge18$ is even and $0\le r\le k$. It suffices to prove that there is an absolute constant $C$ such that $n<Ck$ whenever $G$ has no $K_{k+1}$ component. 

We first establish comparison bounds for $\nu(G)$ and $\Psi(G)$ in Lemma~\ref{lem:linear-comparison}, and then show in Lemma~\ref{lem:periodic-coarse-core} that at most one set of $k+1$ vertices induces a graph with few nonedges.

\begin{lemma}\label{lem:linear-comparison}
If $q\ge2$, $r$ is even and $r>2(k+1)/3$, then
\[
 \nu(G)\le2(k+1)^2(k+1-r)
 \quad \text{and} \quad
 \Psi(G)\ge g_{k+1}(n)-9(k+1)(k+1-r)\quad\text{if }r<k.
\]
\end{lemma}
\begin{proof}
The lower bound in Theorem~\ref{mainbounds} and the identity $\nu(G)=n\binom{k}{2}-3\tau(G)$ from~\eqref{intro_path_identity} give
\[
 \nu(G)\le \frac{r(k+1-r)(k+1+r)}2+\frac{3r(2r-k-2)}4
 \le (k+1)^2(k+1-r)+\frac34(k+1)^2\le2(k+1)^2(k+1-r).
\]
For the second assertion, let $r<k$. We first construct a $k$-regular graph $H$ on $2k+r+2$ vertices as follows. Take three disjoint cliques, one of order $k$ and two of order $(k+r+2)/2$. Using Construction~\ref{construction:regular-bipartite}, join the latter two cliques by a $((k-r)/2)$-regular bipartite graph containing a perfect matching. Choose $k/2$ edges of this perfect matching, delete them, and add a perfect matching between their $k$ endpoints and the vertices of the first clique. Note that, before deletion, the endpoints of each selected matching edge have $(k-r)/2-1$ common neighbors in each of the latter two cliques, giving $k-r-2$ triangles per selected edge, with no triangle counted twice. Moreover, no triangle in $H$ contains an edge of the newly added matching. Therefore
\[
 \begin{aligned}
 \tau(H)=\binom{k}{3}+2\binom{(k+r+2)/2}{3}
 +(k+r+2)\binom{(k-r)/2}{2}-\frac{k(k-r-2)}2
 =2\binom{k}{3}+\binom{r+2}{3}.
 \end{aligned}
\]
Comparing $G$ with $(q-2)K_{k+1}\cup H$ and using $\Psi(K_{k+1})=(k+1)^3$ together with the identity $\Psi(H)=(2k+r+2)(3k+1)+6\tau(H)$ from~\eqref{psireg}, we obtain
\[
 \Psi(G)\ge\Psi\bigl((q-2)K_{k+1}\cup H\bigr)
 =(q-2)\Psi(K_{k+1})+\Psi(H)
 \ge g_{k+1}(n)-9(k+1)(k+1-r).\qedhere
\]
\end{proof}

\begin{lemma}\label{lem:periodic-coarse-core}
Fix an absolute constant $\eta\in(0,1/16)$. If $n>8(k+1)$ and $G$ has no $K_{k+1}$ component, then at most one set of $k+1$ vertices induces a graph with at most $\eta k$ nonedges. Consequently, $\nu(G)>\eta k(n-k-1)$.
\end{lemma}
\begin{proof}
Suppose that there exist two distinct sets $V_1,V_2$ of order $k+1$ satisfy $e(H_i)\le\eta k$, where $H_i=\overline{G[V_i]}$ for $i\in[2]$. Since $G$ has no $K_{k+1}$ component, $e(H_i)>0$. Since $G$ is $k$-regular and $|V_{i}|=k+1$, we have
\begin{equation}\label{eq:linear-external-degree}
 |N_G(x)\setminus V_i|=d_{H_i}(x)\qquad \text{for every } x\in V_i.
\end{equation}

We first claim that $V_1$ and $V_2$ are disjoint. Otherwise $V_1\cap V_2\ne\emptyset$. Since $V_1$ and $V_2$ are distinct, we have $1\le |V_1\cap V_2|\le k$. Since $G$ is $k$-regular, for every $x\in V_1\cap V_2$ the number of nonedges incident with $x$ in $G[V_1\cup V_2]$ satisfies $\bigl|(V_1\cup V_2)\setminus N_G[x]\bigr|\ge |V_1\cup V_2|-(k+1)=k+1-|V_1\cap V_2|$. Summing over $x\in V_1\cap V_2$, we obtain $e(H_1)+e(H_2)\ge\sum_{x\in V_1\cap V_2}\bigl|(V_1\cup V_2)\setminus N_G[x]\bigr|\ge |V_1\cap V_2|(k+1-|V_1\cap V_2|)\ge k$, contradicting $e(H_1)+e(H_2)\le2\eta k<k$. 

Next we claim that no edge joins $V_1$ and $V_2$. Otherwise, choose $v_1v_2\in E(G)$ with $v_i\in V_i$ for $i\in[2]$, and choose $u_i\in V_i\setminus\{v_i\}$ with $v_iu_i\notin E(G)$. For each $i\in[2]$, if $z\in V_i\setminus\{v_i,u_i\}$ is not a common neighbor of $v_i$ and $u_i$ in $G$, then $v_iz\notin E(G)$ or $u_iz\notin E(G)$. Different vertices $z$ yield distinct nonedges of $G[V_i]$ other than $v_iu_i$, so at most $e(H_i)-1$ of the $k-1$ vertices in $V_i\setminus\{v_i,u_i\}$ are not common neighbors of $v_i$ and $u_i$. Hence, adding $v_iu_i$ creates at least $(k-1)-(e(H_i)-1)=k-e(H_i)$ triangles inside $V_i$.

If $u_1u_2\in E(G)$, then we replace $v_1v_2,u_1u_2$ by $v_1u_1,v_2u_2$. Otherwise, since $v_iu_i$ is a nonedge of $G[V_i]$, \eqref{eq:linear-external-degree} gives $|N_G(u_i)\setminus V_i|=d_{H_i}(u_i)\ge1$ for each $i\in[2]$, so we can choose $w_i\in N_G(u_i)\setminus V_i$. Since $G$ is $k$-regular and each of $V_1,V_2,N_G[w_1],N_G[w_2]$ has $k+1$ vertices, at most $4k(k+1)$ edges of $G$ have an endpoint in $V_1\cup V_2\cup N_G[w_1]\cup N_G[w_2]$. As $e(G)=nk/2>4k(k+1)$, we can choose $xy\in E(G)$ with both endpoints outside this union and replace $v_1v_2,u_1w_1,u_2w_2,xy$ by $v_1u_1,v_2u_2,w_1x,w_2y$. Note that both replacements give simple $k$-regular graphs, even when $w_1=w_2$, and delete no edge of $G[V_1]$ or $G[V_2]$.

By~\eqref{eq:linear-external-degree} and our contradiction assumption $e(H_i)\le\eta k$, $e_G(V_i,V(G)\setminus V_i)=2e(H_i)\le2\eta k$ for each $i\in[2]$. Since distinct triangles containing a fixed edge with one endpoint in $V_i$ and the other in $V(G)\setminus V_i$ have distinct second edges between $V_i$ and $V(G)\setminus V_i$, the fixed edge lies in at most $2\eta k$ triangles. The first replacement therefore loses at most $4\eta k$ triangles. Since $xy$ lies in at most $k-1$ triangles, the second replacement loses at most $6\eta k+(k-1)$ triangles. Thus the loss is at most $6\eta k+(k-1)$ in either case, while adding $v_1u_1,v_2u_2$ creates at least $(k-e(H_1))+(k-e(H_2))\ge2(1-\eta)k$ triangles. Hence the net increase is at least $2(1-\eta)k-6\eta k-(k-1)=(1-8\eta)k+1>0$, contradicting the extremality of $G$. 

We now construct a graph $G'$ on $V(G)$. Let $R=V(G)\setminus(V_1\cup V_2)$. Since the edges of $H_1,H_2$ have at most $2(e(H_1)+e(H_2))\le4\eta k<k+1$ distinct endpoints, we can choose $W\subseteq V_1\cup V_2$ of order $k+1$ containing all these endpoints. By~\eqref{eq:linear-external-degree}, every vertex of $V_1\cup V_2$ with a neighbor in $R$ lies in $W$. Let $G'$ be the graph obtained from $G$ by making $(V_1\cup V_2)\setminus W$ induce $K_{k+1}$, taking $E(H_1)\cup E(H_2)$ as the nonedges in $W$, and deleting all edges between $W$ and $(V_1\cup V_2)\setminus W$, while leaving all edges with at least one endpoint in $R$ unchanged. Equation~\eqref{eq:linear-external-degree} shows that $G'$ is $k$-regular.

For any graph $H$, inclusion--exclusion gives
\begin{equation}\label{eq:linear-complement-triangles}
 \tau(\overline H)=\binom{|V(H)|}{3}-(|V(H)|-2)e(H)
       +\sum_{x\in V(H)}\binom{d_H(x)}2-\tau(H).
\end{equation}
Applying this identity to $H_1$, $H_2$, and $\overline{G'[W]}$ gives $\tau(G'[W])=\tau(G[V_1])+\tau(G[V_2])-\binom{k+1}{3}$. Adding the $\binom{k+1}{3}$ triangles in $G'[(V_1\cup V_2)\setminus W]$ then gives $\tau\bigl(G'[V_1\cup V_2]\bigr)=\tau\bigl(G[V_1\cup V_2]\bigr)$. Note that every triangle of $G$ meeting $R$ has all its vertices outside $R$ in $W$. Since every edge of $G[W]$ is present in $G'[W]$ and every edge of $G$ with at least one endpoint in $R$ is retained in $G'$, every triangle of $G$ meeting $R$ is also a triangle of $G'$. Extremality therefore gives $\tau(G')=\tau(G)$, so $G'$ is also extremal.

We claim that no vertex of $R$ has neighbors in both $V_1$ and $V_2$ in $G$. Otherwise, choose a vertex $v\in R$ adjacent in $G$ to both $x_1\in V_1$ and $x_2\in V_2$. By~\eqref{eq:linear-external-degree}, we have $x_1,x_2\in W$. The edge $x_1x_2$ added in $G'$ then creates a triangle $vx_1x_2$ that is absent from $G$, contradicting $\tau(G')=\tau(G)$.

We now construct a $k$-regular graph $G''$ on $V(G)$. For each $i\in[2]$, choose $z_i\in V_i$ with $d_{H_i}(z_i)>0$. Let $G''$ be the graph obtained from $G'$ by replacing $z_1x$ with $z_2x$ for every $x\in N_{H_2}(z_2)$ and $z_2y$ with $z_1y$ for every $y\in N_{G'}(z_2)\cap R$. Since $|N_{G'}(z_2)\cap R|=d_{H_2}(z_2)$ by~\eqref{eq:linear-external-degree}, these replacements preserve every vertex degree. Note that all added edges are absent from $G'$, so $G''$ is a simple $k$-regular graph.

The clique on $(V_1\cup V_2)\setminus W$ is unchanged when constructing $G''$ from $G'$. For distinct $u,v\in W\setminus\{z_1,z_2\}$, the number of triangles among $z_1uv,z_2uv$ is the same in $G'$ and $G''$, except when one of $u,v$ lies in $N_{H_1}(z_1)$ and the other in $N_{H_2}(z_2)$, in which case it increases from zero to one. Triangles inside $W$ containing both $z_1,z_2$ or neither are unchanged. Hence $\tau(G''[W])-\tau(G'[W])=d_{H_1}(z_1)d_{H_2}(z_2)$.

Since no vertex of $R$ has neighbors in both $V_1$ and $V_2$ in $G$, triangles of $G'$ containing $z_1$ and meeting $R$ are retained in $G''$, while those containing $z_2$ and meeting $R$ become triangles of $G''$ upon replacing $z_2$ by $z_1$. The latter triangles are distinct and were absent from $G'$, since $z_1,z_2$ have no common neighbor in $R$ in $G'$. All other triangles of $G'$ meeting $R$ are unchanged. Thus the number of triangles meeting $R$ does not decrease, and $\tau(G'')-\tau(G')\ge d_{H_1}(z_1)d_{H_2}(z_2)>0$. This contradicts the extremality of $G'$ and proves the first assertion of Lemma~\ref{lem:periodic-coarse-core}.

It holds that for at most $k+1$ vertices $v\in V(G)$, the subgraph of $G$ induced by $v$ and its neighbors has at most $\eta k$ nonedges. Every other vertex is the center of more than $\eta k$ induced copies of $P_3$. Summing over their centers yields $\nu(G)>\eta k(n-k-1)$, as deirsed.
\end{proof}

\begin{proof}[\textbf{Proof of Theorem~\ref{mainperiodicity}}]
Fix $\eta\in(0,1/16)$ as in Lemma~\ref{lem:periodic-coarse-core}. Choose an integer $M>1+9/\eta$, and then fix $0<\varepsilon<\min\{1/100,1/(20M)\}$.
All these constants are independent of $k,n,r$.

Assume that $G$ has no $K_{k+1}$ component. Since $q<2M$ gives $n<2M(k+1)=O(k)$, we may assume $q\ge2M$, so $n>8(k+1)$. By Theorem~\ref{mainodd}, we may also assume that $r>2(k+1)/3$ is even. If $k+1-r\le\varepsilon^{-1}$, Lemmas~\ref{lem:linear-comparison} and~\ref{lem:periodic-coarse-core} give $n=O(k)$. Thus we may assume $k+1-r>\varepsilon^{-1}$, so $r<k$, and Lemma~\ref{lem:linear-comparison} gives
\begin{equation}\label{eq:short-periodic-comparison}
 \nu(G)\le 2(k+1)^2(k+1-r),\qquad
 \Psi(G)\ge g_{k+1}(n)-9(k+1)(k+1-r).
\end{equation}

Choose a maximal family $\{U_i\}$ of pairwise disjoint subsets of $V(G)$ satisfying
\[
 k+1-\varepsilon(k+1-r)\le |U_i|\le k \quad \text{and} \quad
 \Delta(\overline{G[U_i]})\le\varepsilon(k+1),
\]
and, for each member $U_i$ of this family, let
\[
 B_i=\{v\in V(G):|N_G(v)\cap U_i|>2k/3\}.
\]
\begin{observation}\label{obs:short-periodic-basic}
For each $i$, the following hold. \begin{enumerate*}[label=$(\roman*)$,font=\normalfont,ref=\thetheorem(\roman*),mode=unboxed,itemjoin=\space,labelsep=0.4em]\item\label{obs:short-periodic-basic-i} $U_i\subseteq B_i$; \item\label{obs:short-periodic-basic-ii} $B_i\cap B_j=\emptyset$ for every $j\ne i$; \item\label{obs:short-periodic-basic-iii} $|B_i|-|U_i|<2\varepsilon(k+1)$; \item\label{obs:short-periodic-basic-iv} $|B_i|<(1+2\varepsilon)(k+1)$; \item\label{obs:short-periodic-basic-v} $d_{\overline{G[B_i]}}(x)\le3\varepsilon(k+1)$ for every $x\in U_i$; \item\label{obs:short-periodic-basic-vi} $|N_G(x)\setminus B_i|<k/3$ for every $x\in B_i\setminus U_i$. \end{enumerate*}
\end{observation}
\begin{proof}
For $x\in U_i$, $|N_G(x)\setminus U_i|\le k+1-|U_i|+\varepsilon(k+1)\le\varepsilon(2k+2-r)<k/3$, proving $(i)$. Since $U_i\cap U_j=\emptyset$ for $i\ne j$, any $x\in B_i\cap B_j$ would satisfy $k=d_G(x)\ge |N_G(x)\cap U_i|+|N_G(x)\cap U_j|>4k/3$, a contradiction. This proves $(ii)$. Each vertex of $B_i\setminus U_i$ has more than $2k/3$ neighbors in $U_i$, while each vertex of $U_i$ has at most $\varepsilon(2k+2-r)$ neighbors outside $U_i$. Counting edges between $U_i$ and $B_i\setminus U_i$ therefore gives $(2k/3)(|B_i|-|U_i|)\le e_G(U_i,B_i\setminus U_i)\le |U_i|\varepsilon(2k+2-r)<(4k/3)\varepsilon(k+1)$, where the last inequality uses $|U_i|\le k$ and $r>2(k+1)/3$. This proves $(iii)$. Now $(iv)$ follows from $(iii)$ and $|U_i|\le k$. For $x\in U_i$, at most $\varepsilon(k+1)$ vertices in $U_i\setminus\{x\}$ and at most $|B_i|-|U_i|$ vertices in $B_i\setminus U_i$ are nonadjacent to $x$ in $G$. Hence $(iii)$ gives $d_{\overline{G[B_i]}}(x)\le\varepsilon(k+1)+|B_i|-|U_i|<3\varepsilon(k+1)$, proving $(v)$. Finally, $(vi)$ follows immediately from $(i)$ and the definition of $B_i$.
\end{proof}

For each $i$ and $x\in B_i$, let $d_i(x)=|N_G(x)\setminus B_i|$. Since $G$ is $k$-regular, we have
\begin{equation}\label{eq:short-periodic-external-degree}
 d_i(x)=k+1-|B_i|+d_{\overline{G[B_i]}}(x).
\end{equation}
We first establish the counting inequality in Claim~\ref{clm:short-periodic-count} for use in the proofs of Claims~\ref{clm:short-periodic-large} and~\ref{clm:short-periodic-small}.

\begin{claim}\label{clm:short-periodic-count}
For any $m$ selected sets $B_i$, write them as $B_1,\ldots,B_m$ and relabel the corresponding $U_i,d_i$ accordingly. Let $R=V(G)\setminus\bigcup_{i=1}^mB_i$. Then a lower bound for $\Psi(G[R])+m(k+1)^3-\Psi(G)$ is
\begin{align*}
 \sum_{i=1}^m(k+1-|B_i|)^3
  +3\sum_{i=1}^m\sum_{x\in B_i}d_i(x)\bigl(|B_i|-2d_{\overline{G[B_i]}}(x)\bigr)
-3\sum_{i=1}^m\sum_{v\in R}|N_G(v)\cap B_i|^2.
\end{align*}
\end{claim}
\begin{proof}
Applying the complement identity~\eqref{eq:linear-complement-triangles} to $\overline{G[B_i]}$ gives
\[
 \tau(G[B_i])=\binom{|B_i|}{3}-(|B_i|-2)e(\overline{G[B_i]})+\sum_{x\in B_i}\binom{d_{\overline{G[B_i]}}(x)}{2}-\tau(\overline{G[B_i]}).
\]
Using the preceding identity together with the definition of $\Psi$ and the handshaking lemma, we obtain
\[
 \begin{aligned}
 \Psi(G[B_i])
 &=|B_i|+6e(G[B_i])+6\tau(G[B_i])\\
 &=|B_i|+6\binom{|B_i|}{2}-6e(\overline{G[B_i]})+6\binom{|B_i|}{3}-6(|B_i|-2)e(\overline{G[B_i]})
+6\sum_{x\in B_i}\binom{d_{\overline{G[B_i]}}(x)}{2}-6\tau(\overline{G[B_i]})\\
 &=|B_i|^3-3(|B_i|-1)\sum_{x\in B_i}d_{\overline{G[B_i]}}(x)
+3\sum_{x\in B_i}d_{\overline{G[B_i]}}(x)\bigl(d_{\overline{G[B_i]}}(x)-1\bigr)-6\tau(\overline{G[B_i]})\\
 &=|B_i|^3-3\sum_{x\in B_i}d_{\overline{G[B_i]}}(x)\bigl(|B_i|-d_{\overline{G[B_i]}}(x)\bigr)-6\tau(\overline{G[B_i]}).
 \end{aligned}
\]
By Observation~\ref{obs:short-periodic-basic-ii}, $B_i\cap B_j=\emptyset$ for distinct $i,j\in[m]$. Recall that $\Psi(G)$ counts ordered triples of vertices, with repetitions allowed, in which any two distinct vertices are adjacent. Every such triple not contained in one of $B_1,\ldots,B_m,R$ contains a vertex $x\in B_i$ with the other two vertices outside $B_i$, or a vertex $v\in R$ with the other two vertices in the same $B_i$. Choosing the other two vertices in order and allowing $x$ or $v$ to occur first, second, or third gives at most $3d_i(x)^2$ triples for each $x\in B_i$ in the first case and $3|N_G(v)\cap B_i|^2$ triples for each $v\in R$ in the second. Hence
\[
 \begin{aligned}
 \Psi(G)&\le\Psi(G[R])+\sum_{i=1}^m\Psi(G[B_i])
 +3\sum_{i=1}^m\sum_{x\in B_i}d_i(x)^2
 +3\sum_{i=1}^m\sum_{v\in R}|N_G(v)\cap B_i|^2.
 \end{aligned}
\]
For each $i$, substituting~\eqref{eq:short-periodic-external-degree} into the expression for $\Psi(G[B_i])$ above and simplifying gives
\[
 \begin{aligned}
 (k+1)^3-\Psi(G[B_i])-3\sum_{x\in B_i}d_i(x)^2
=(k+1-|B_i|)^3+3\sum_{x\in B_i}d_i(x)\bigl(|B_i|-2d_{\overline{G[B_i]}}(x)\bigr)+6\tau(\overline{G[B_i]}).
 \end{aligned}
\]
Summing over $i$, using the preceding bound on $\Psi(G)$, and discarding the nonnegative terms $6\tau(\overline{G[B_i]})$ proves the claim.
\end{proof}

\begin{claim}\label{clm:short-periodic-large}
At most two sets $B_i$ have size greater than $k+1$.
\end{claim}
\begin{proof}
Suppose that, after relabeling the indices, $B_1,B_2,B_3$ all have size greater than $k+1$. Let $R=V(G)\setminus\bigcup_{i=1}^3 B_i$ and $s=\sum_{i=1}^3(|B_i|-k-1)$. Clearly, $s\ge3$.

For each $i\in[3]$, Observation~\ref{obs:short-periodic-basic-i} gives $U_i\subseteq B_i$, and Observation~\ref{obs:short-periodic-basic-iii} gives $|B_i|-|U_i|<2\varepsilon(k+1)<k/24$. Thus, for $v\in R$, the definition of $B_i$ gives $|N_G(v)\cap U_i|\le2k/3$ and $|N_G(v)\cap(B_i\setminus U_i)|<k/24$. Using these neighbor bounds we obtain, for $v\in R$,
\[
 \begin{aligned}
 |N_G(v)\cap B_i|^2
 &=|N_G(v)\cap U_i|^2+2|N_G(v)\cap U_i|\,|N_G(v)\cap(B_i\setminus U_i)|+|N_G(v)\cap(B_i\setminus U_i)|^2\\
 &\le\frac{3k}{4}|N_G(v)\cap U_i|+\frac{k}{4}|N_G(v)\cap(B_i\setminus U_i)|.
 \end{aligned}
\]
For $x\in U_i$, Observation~\ref{obs:short-periodic-basic-v} gives $d_{\overline{G[B_i]}}(x)\le3\varepsilon(k+1)$, so $|B_i|-2d_{\overline{G[B_i]}}(x)>3k/4$. For $x\in B_i\setminus U_i$, Observation~\ref{obs:short-periodic-basic-vi} gives $d_i(x)<k/3$. Together with $|B_i|<(1+2\varepsilon)(k+1)$ from Observation~\ref{obs:short-periodic-basic-iv}, identity~\eqref{eq:short-periodic-external-degree} gives $|B_i|-2d_{\overline{G[B_i]}}(x)=2(k+1)-|B_i|-2d_i(x)>k/4$. Summing the preceding square bound over $v\in R$, counting edges, and using $R\subseteq V(G)\setminus B_i$ and these coefficient bounds, we obtain
\[
 \begin{aligned}
 \sum_{v\in R}|N_G(v)\cap B_i|^2
 &\le\frac{3k}{4}e_G(U_i,R)+\frac{k}{4}e_G(B_i\setminus U_i,R)\\
 &\le\frac{3k}{4}\sum_{x\in U_i}d_i(x)+\frac{k}{4}\sum_{x\in B_i\setminus U_i}d_i(x)\\
 &\le\sum_{x\in B_i}d_i(x)\bigl(|B_i|-2d_{\overline{G[B_i]}}(x)\bigr).
 \end{aligned}
\]

Using this inequality in Claim~\ref{clm:short-periodic-count} with $m=3$, and using $\sum_{i=1}^3(|B_i|-k-1)^3\le s^3$, we obtain
\[
 \Psi(G)\le\Psi(G[R])+3(k+1)^3+\sum_{i=1}^3(|B_i|-k-1)^3\le\Psi(G[R])+3(k+1)^3+s^3.
\]

By Observation~\ref{obs:short-periodic-basic-ii}, $B_i\cap B_j=\emptyset$ for $i\ne j$, so $|R|=(q-3)(k+1)+(r-s)$. Observation~\ref{obs:short-periodic-basic-iv} gives $|B_i|<(1+2\varepsilon)(k+1)$ for each $i\in[3]$. Since $\varepsilon<1/100$, we have $s<6\varepsilon(k+1)<(k+1)/10$, which together with $r>2(k+1)/3$ gives $(k+1)/2<r-s<k+1$. Since $q\ge2M>3$ and $G$ is $k$-regular, we have $q-3\ge0$ and $\Delta(G[R])\le k$. Thus Lemma~\ref{psilemma} applies to $G[R]$ and gives $\Psi(G[R])\le(q-3)(k+1)^3+(r-s)^3$. Substituting this into the preceding bound on $\Psi(G)$ yields
\[
 \Psi(G)\le q(k+1)^3+(r-s)^3+s^3=g_{k+1}(n)-3sr(r-s).
\]

Since $r>2(k+1)/3$, $r-s>(k+1)/2$, and $s\ge3$, we have $3sr(r-s)>s(k+1)^2\ge3(k+1)^2>9(k+1)(k+1-r)$, contradicting~\eqref{eq:short-periodic-comparison}. Thus at most two sets $B_i$ have size greater than $k+1$.
\end{proof}

\begin{claim}\label{clm:short-periodic-small}
Fewer than $M$ sets $B_i$ have size at most $k+1$.
\end{claim}
\begin{proof}
Suppose that $B_1,\ldots,B_M$ all have size at most $k+1$, after relabeling. Let $R=V(G)\setminus\bigcup_{i=1}^M B_i$, $t=\sum_{i=1}^M(k+1-|B_i|)$, 
\[
 D=\sum_{i=1}^M\sum_{x\in U_i}d_{\overline{G[B_i]}}(x),\quad\text{and}\quad
 Q=\sum_{i=1}^M\sum_{v\in R}|N_G(v)\cap U_i|^2.
\]
By Observation~\ref{obs:short-periodic-basic-i}, $U_i\subseteq B_i$, so $k+1-|B_i|\le k+1-|U_i|\le\varepsilon(k+1-r)$ and hence $t\le M\varepsilon(k+1-r)<(k+1-r)/4$. We establish lower and upper bounds for $\Psi(G[R])+3Q$.

\noindent\textbf{Lower bound.} Since $U_i\subseteq B_i$ and $|N_G(v)\cap(B_i\setminus U_i)|\le |B_i|-|U_i|$, for $v\in R$
\begin{equation}\label{eq:short-periodic-square-difference}
 \begin{aligned}[b]
 &\sum_{v\in R}\bigl(|N_G(v)\cap B_i|^2-|N_G(v)\cap U_i|^2\bigr)\\
 &=2\sum_{v\in R}|N_G(v)\cap U_i|\,|N_G(v)\cap(B_i\setminus U_i)|+\sum_{v\in R}|N_G(v)\cap(B_i\setminus U_i)|^2\\
 &\le2(|B_i|-|U_i|)\sum_{v\in R}|N_G(v)\cap U_i|+(|B_i|-|U_i|)\sum_{v\in R}|N_G(v)\cap(B_i\setminus U_i)|\\
 &=2(|B_i|-|U_i|)\sum_{x\in U_i}|N_G(x)\cap R|+(|B_i|-|U_i|)\sum_{x\in B_i\setminus U_i}|N_G(x)\cap R|\\
 &\le2(|B_i|-|U_i|)\sum_{x\in U_i}d_i(x)+(|B_i|-|U_i|)\sum_{x\in B_i\setminus U_i}d_i(x).
 \end{aligned}
\end{equation}

Next, $|U_i|\ge k+1-\varepsilon(k+1-r)$ gives $|U_i|^2\ge(k+1)^2-2\varepsilon(k+1)(k+1-r)$. Since $k+1-|B_i|\ge0$, we have
\begin{equation}\label{eq:short-periodic-deficit-bound}
 3(k+1-|B_i|)|U_i|^2\ge\bigl(3(k+1)^2-6\varepsilon(k+1)(k+1-r)\bigr)(k+1-|B_i|).
\end{equation}
For $x\in U_i$, the inclusion $U_i\subseteq B_i$ and the defining conditions on $U_i$ give $d_i(x)\le k+1-|U_i|+\varepsilon(k+1)$. Since $|B_i|\le k+1$, we have $2|U_i|-|B_i|-2d_i(x)\ge4|U_i|-3(k+1)-2\varepsilon(k+1)\ge5(k+1)/6$. For $x\in B_i\setminus U_i$, Observation~\ref{obs:short-periodic-basic-vi} gives $d_i(x)<k/3$, so $|U_i|-2d_i(x)>|U_i|-2k/3\ge(k+1)/6$. Combining these bounds with the handshaking lemma, we obtain
\begin{equation}\label{eq:short-periodic-degree-bound}
 \begin{aligned}[b]
 &3\sum_{x\in U_i}d_{\overline{G[B_i]}}(x)(2|U_i|-|B_i|-2d_i(x))+3\sum_{x\in B_i\setminus U_i}d_{\overline{G[B_i]}}(x)(|U_i|-2d_i(x))\\
 &\ge2(k+1)\sum_{x\in U_i}d_{\overline{G[B_i]}}(x)+(k+1)e(\overline{G[B_i]})\\
 &\ge2k\sum_{x\in U_i}d_{\overline{G[B_i]}}(x)+(k+1)e(\overline{G[B_i]}).
 \end{aligned}
\end{equation}

Finally, apply Claim~\ref{clm:short-periodic-count} with $m=M$ and add $3Q$ to both sides. Since the square-sum difference has coefficient $-3$, the upper bound in~\eqref{eq:short-periodic-square-difference} gives a lower bound for $\Psi(G[R])+3Q+M(k+1)^3-\Psi(G)$. More explicitly, this lower bound is the sum over $i$ of the following expressions, which we rewrite by grouping the terms over $U_i$ and $B_i\setminus U_i$ and using~\eqref{eq:short-periodic-external-degree}.
\[
 \begin{aligned}
 &(k+1-|B_i|)^3+3\sum_{x\in B_i}d_i(x)\bigl(|B_i|-2d_{\overline{G[B_i]}}(x)\bigr)
-6(|B_i|-|U_i|)\sum_{x\in U_i}d_i(x)-3(|B_i|-|U_i|)\sum_{x\in B_i\setminus U_i}d_i(x)\\
 &=(k+1-|B_i|)^3+3\sum_{x\in U_i}d_i(x)\bigl(2|U_i|-|B_i|-2d_{\overline{G[B_i]}}(x)\bigr)
 +3\sum_{x\in B_i\setminus U_i}d_i(x)\bigl(|U_i|-2d_{\overline{G[B_i]}}(x)\bigr)\\
 &=(k+1-|B_i|)^3+3\sum_{x\in U_i}(k+1-|B_i|)(2|U_i|-|B_i|)
 +3\sum_{x\in B_i\setminus U_i}(k+1-|B_i|)|U_i|\\
 &\quad+3\sum_{x\in U_i}d_{\overline{G[B_i]}}(x)(2|U_i|-|B_i|-2d_i(x))+3\sum_{x\in B_i\setminus U_i}d_{\overline{G[B_i]}}(x)(|U_i|-2d_i(x))\\
 &=(k+1-|B_i|)^3+3(k+1-|B_i|)|U_i|^2\\
 &\quad+3\sum_{x\in U_i}d_{\overline{G[B_i]}}(x)(2|U_i|-|B_i|-2d_i(x))
+3\sum_{x\in B_i\setminus U_i}d_{\overline{G[B_i]}}(x)(|U_i|-2d_i(x)).
 \end{aligned}
\]
Summing the last expression above over $i\in[M]$, using~\eqref{eq:short-periodic-deficit-bound} and~\eqref{eq:short-periodic-degree-bound}, and discarding the nonnegative terms $(k+1-|B_i|)^3$, we obtain
\begin{equation}\label{eq:short-periodic-lower}
 \Psi(G[R])+3Q+M(k+1)^3-\Psi(G)
 \ge(3(k+1)^2-6\varepsilon(k+1)(k+1-r))t+2kD+(k+1)\sum_{i=1}^M e(\overline{G[B_i]}).
\end{equation}

\noindent\textbf{Upper bound.} We first estimate $Q$. Counting edges from the pairwise disjoint sets $U_i$ to $R\subseteq V(G)\setminus B_i$ and using~\eqref{eq:short-periodic-external-degree}, the bounds $|U_i|\le k$ from the choice of $U_i$ and $|B_i|\le k+1$ from our supposition, and the definitions of $D,t$, we obtain
\begin{equation}\label{eq:short-periodic-edge-bound}
 \begin{aligned}[b]
 e_G\Bigl(R,\bigcup_{i=1}^M U_i\Bigr)
 =\sum_{i=1}^M\sum_{x\in U_i}|N_G(x)\cap R|
 \le\sum_{i=1}^M\sum_{x\in U_i}d_i(x)
 =\sum_{i=1}^M|U_i|(k+1-|B_i|)+D\le kt+D.
 \end{aligned}
\end{equation}

For $v\in R$ and $i\in[M]$, the definition of $B_i$ and the fact that $G$ is $k$-regular give $|N_G(v)\cap U_i|\le\min\{2k/3,k-d_{G[R]}(v)\}$. Since the sets $U_i$ are pairwise disjoint, it follows that
\begin{equation}\label{eq:short-periodic-q-sum}
 Q=\sum_{i=1}^M\sum_{v\in R}|N_G(v)\cap U_i|^2\le\sum_{v\in R}\min\{2k/3,k-d_{G[R]}(v)\}\left|N_G(v)\cap\bigcup_{i=1}^M U_i\right|.
\end{equation}

By Observation~\ref{obs:short-periodic-basic-ii}, the sets $B_i$ are pairwise disjoint, so $|R|=(q-M)(k+1)+r+t\ge t$. Choose $t$ vertices of smallest degree in $G[R]$, and let $W$ and $W'$ consist of those with degree at most $r-1$ and at least $r$, respectively. Then $|W|+|W'|=t$.

Since $G$ is $k$-regular, $\left|N_G(v)\cap\bigcup_{i=1}^M U_i\right|\le k$ for every $v\in R$, and~\eqref{eq:short-periodic-edge-bound} bounds the sum of these counts over $v\in R$ by $kt+D$. The $t$ largest values of $\min\{2k/3,k-d_{G[R]}(v)\}$ occur on $W\cup W'$, and all values are at most $2k/3$. Applying these bounds to~\eqref{eq:short-periodic-q-sum} gives
\[
 Q\le k\sum_{v\in W\cup W'}\min\{2k/3,k-d_{G[R]}(v)\}+\frac{2k}{3}D.
\]
For $v\in W'$, the inequality $d_{G[R]}(v)\ge r$ gives $\min\{2k/3,k-d_{G[R]}(v)\}\le k-r<k+1-r$. Splitting the sum over $W\cup W'$ therefore yields
\begin{equation}\label{eq:short-periodic-q-bound}
 Q\le k\sum_{v\in W}\min\{2k/3,k-d_{G[R]}(v)\}+k(k+1-r)|W'|+\frac{2k}{3}D.
\end{equation}

Next, deleting $v\in W$ from $G[R]$ removes $v$, its $d_{G[R]}(v)$ incident edges, and at most $\binom{d_{G[R]}(v)}2$ triangles containing $v$. By the definition of $\Psi$, each deleted vertex contributes $1$ and each deleted edge or triangle $6$. Summing over $v\in W$ counts every deleted edge and triangle at least once, giving
\begin{equation}\label{eq:short-periodic-deletion-bound}
 \begin{aligned}[b]
 \Psi(G[R])-\Psi(G[R\setminus W])
 \le\sum_{v\in W}\left(1+6d_{G[R]}(v)+6\binom{d_{G[R]}(v)}2\right)
 =\sum_{v\in W}\bigl(1+3d_{G[R]}(v)(d_{G[R]}(v)+1)\bigr).
 \end{aligned}
\end{equation}

We add the upper bound~\eqref{eq:short-periodic-deletion-bound} on the loss of $\Psi$ from deleting $W$ to three times the sum over $v\in W$ in~\eqref{eq:short-periodic-q-bound}. Since $d_{G[R]}(v)\le r-1$, we have $d_{G[R]}(v)(d_{G[R]}(v)+1)\le r\,d_{G[R]}(v)$. It holds that
\begin{equation}\label{eq:short-periodic-scalar}
 \begin{aligned}[b]
 &1+3d_{G[R]}(v)(d_{G[R]}(v)+1)+3k\min\{2k/3,k-d_{G[R]}(v)\}\\
 &\le1+3r\,d_{G[R]}(v)+3r\min\{2k/3,k-d_{G[R]}(v)\}+3(k-r)\min\{2k/3,k-d_{G[R]}(v)\}\\
 &\le1+3r\,d_{G[R]}(v)+3r(k-d_{G[R]}(v))+2k(k-r)\\
 &<3(k+1)^2-(k+1)(k+1-r).
 \end{aligned}
\end{equation}

Combining~\eqref{eq:short-periodic-deletion-bound} with three times~\eqref{eq:short-periodic-q-bound} and applying~\eqref{eq:short-periodic-scalar} for each $v\in W$, we obtain
\begin{equation}\label{eq:short-periodic-combined-bound}
 \Psi(G[R])+3Q\le\Psi(G[R\setminus W])+\bigl(3(k+1)^2-(k+1)(k+1-r)\bigr)|W|+3k(k+1-r)|W'|+2kD.
\end{equation}
We next bound $\Psi(G[R\setminus W])+3k(k+1-r)|W'|$. The graph $G[R\setminus W]$ has $(q-M)(k+1)+r+|W'|$ vertices and maximum degree at most $k$, with $q-M\ge0$ and $0\le r+|W'|\le r+t<k+1$, so Lemma~\ref{psilemma} gives $\Psi(G[R\setminus W])\le(q-M)(k+1)^3+(r+|W'|)^3$. Since $0\le |W'|\le t<(k+1-r)/4<r$, we have $3r|W'|+|W'|^2\le4r|W'|\le r(k+1-r)$. Expanding the cube therefore gives $\Psi(G[R\setminus W])\le(q-M)(k+1)^3+r^3+r(2r+k+1)|W'|$. For the second term, $2(k+1)/3<r\le k$ gives $3k<2(k+1+r)$, so $3k(k+1-r)|W'|\le2(k+1+r)(k+1-r)|W'|$. Adding these two bounds we obtain
\begin{equation}\label{eq:short-periodic-remainder-bound}
 \Psi(G[R\setminus W])+3k(k+1-r)|W'|\le(q-M)(k+1)^3+r^3+\bigl(3(k+1)^2-(k+1)(k+1-r)\bigr)|W'|.
\end{equation}
Substituting~\eqref{eq:short-periodic-remainder-bound} into~\eqref{eq:short-periodic-combined-bound} and using $|W|+|W'|=t$ gives
\begin{equation}\label{eq:short-periodic-upper}
 \begin{aligned}[b]
 \Psi(G[R])+3Q
\leq(q-M)(k+1)^3+r^3+\bigl(3(k+1)^2-(k+1)(k+1-r)\bigr)t+2kD.
 \end{aligned}
\end{equation}

Finally, substituting~\eqref{eq:short-periodic-upper} into~\eqref{eq:short-periodic-lower} cancels $2kD$. Dividing by $k+1$ and using~\eqref{eq:short-periodic-comparison}, we obtain
\[
 (1-6\varepsilon)(k+1-r)\sum_{i=1}^M(k+1-|B_i|)+\sum_{i=1}^M e(\overline{G[B_i]}) \le\frac{g_{k+1}(n)-\Psi(G)}{k+1}\le9(k+1-r).
\]
However, if $|B_i|\le k$, the inequalities $k+1-|B_i|\ge1$ and $1-6\varepsilon>\eta$ imply that $(1-6\varepsilon)(k+1-r)(k+1-|B_i|)>\eta(k+1-r)$. If $|B_i|=k+1$, then $k+1-|B_i|=0$, and since the sets $B_i$ are pairwise disjoint by Observation~\ref{obs:short-periodic-basic-ii}, Lemma~\ref{lem:periodic-coarse-core} gives the nonedge bound $e(\overline{G[B_i]})>\eta k>\eta(k+1-r)$ for all but at most one such index. Since at least $M-1$ indices contribute more than $\eta(k+1-r)$, all terms are nonnegative, and $M>1+9/\eta$, the left-hand side exceeds $(M-1)\eta(k+1-r)>9(k+1-r)$, a contradiction.
\end{proof}

\begin{claim}\label{clm:short-periodic-residual}
Let $R=V(G)\setminus\bigcup_i B_i$, where the union is over all sets $B_i$. Then $e(\overline{G[N_G[v]]})>\varepsilon^2(k+1)(k+1-r)/2$ for every $v\in R$. Consequently, $|R|\le4\varepsilon^{-2}(k+1)$.
\end{claim}
\begin{proof}
For each $v\in V(G)$, let $F_v=\overline{G[N_G[v]]}$. Fix $v\in R$, and let $X=\{x\in N_G[v]:d_{F_v}(x)\le\varepsilon(k+1)\}$. We first show that $X$ is disjoint from every $U_i$. Suppose that $x\in X\cap U_i$. By Observation~\ref{obs:short-periodic-basic-i}, $U_i\subseteq B_i$, so $v\notin U_i$ and $x\in N_G(v)$. Since $G$ is $k$-regular, $d_{F_v}(x)=|N_G(x)\setminus N_G[v]|$. Splitting $U_i\setminus N_G(v)$ according to adjacency to $x$ in $G$ gives
\[
 |U_i\setminus N_G(v)|\le d_{\overline{G[U_i]}}(x)+|N_G(x)\setminus N_G[v]|\le\varepsilon(k+1)+d_{F_v}(x).
\]
Using the preceding estimate and the defining conditions on $U_i$ and $X$, we obtain
\[
 \begin{aligned}
 |N_G(v)\cap U_i|=|U_i|-|U_i\setminus N_G(v)|\ge |U_i|-\varepsilon(k+1)-d_{F_v}(x)
 \ge k+1-\varepsilon(k+1-r)-2\varepsilon(k+1)>2k/3,
 \end{aligned}
\]
contradicting $v\notin B_i$.

Since $X$ is disjoint from every $U_i$, $\Delta(\overline{G[X]})\le\varepsilon(k+1)$, and $\varepsilon(k+1-r)>1$, the maximality of the family $\{U_i\}$ gives $|X|<k+1-\varepsilon(k+1-r)$. Otherwise, a subset of $X$ of size $\min\{|X|,k\}$ could be added to the family.

Every vertex of $N_G[v]\setminus X$ has degree greater than $\varepsilon(k+1)$ in $F_v$. Since $|N_G[v]|=k+1$, the handshaking lemma and the bound $|X|<k+1-\varepsilon(k+1-r)$ give
\begin{equation}\label{eq:short-periodic-residual-nonedges}
 \begin{aligned}[b]
 2e(F_v)\ge\sum_{x\in N_G[v]\setminus X}d_{F_v}(x) >\varepsilon(k+1)|N_G[v]\setminus X|=\varepsilon(k+1)(k+1-|X|)>\varepsilon^2(k+1)(k+1-r).
 \end{aligned}
\end{equation}

For each $v\in V(G)$, the edges of $F_v$ correspond to the induced copies of $P_3$ centered at $v$, so
\[
 \sum_{v\in R}e(F_v)\le\sum_{v\in V(G)}e(F_v)=\nu(G).
\]
Combining this with~\eqref{eq:short-periodic-residual-nonedges} and~\eqref{eq:short-periodic-comparison} gives $(1/2)|R|\varepsilon^2(k+1)(k+1-r)\le\nu(G)\le2(k+1)^2(k+1-r)$, and hence $|R|\le4\varepsilon^{-2}(k+1)$.
\end{proof}

Claims~\ref{clm:short-periodic-large} and~\ref{clm:short-periodic-small} show that fewer than $2M$ sets $B_i$ occur. Observation~\ref{obs:short-periodic-basic-iv} gives $|B_i|<(1+2\varepsilon)(k+1)$ for each $i$. Together with Claim~\ref{clm:short-periodic-residual}, these bounds yield
\[
 n<2M(1+2\varepsilon)(k+1)+4\varepsilon^{-2}(k+1)=O(k).
\]
If $q<2M$, then $n<2M(k+1)=O(k)$. For $q\ge2M$, Theorem~\ref{mainodd} forces a $K_{k+1}$ component if $r$ is odd or $r\le2(k+1)/3$. In the remaining range, Lemmas~\ref{lem:linear-comparison} and~\ref{lem:periodic-coarse-core} give $n=O(k)$ when $k+1-r\le\varepsilon^{-1}$. Together with the bound above, these reductions show that, for some absolute $C$, every extremal $k$-regular graph on at least $Ck$ vertices contains a $K_{k+1}$ component. This completes the proof of the assertion stated at the start of this section. 

Note that for $n\ge Ck$, deleting a $K_{k+1}$ component from an extremal $k$-regular graph on $n+k+1$ vertices leaves an extremal graph on $n$ vertices. The recurrence and the characterization of extremal graphs follow.
\end{proof}
\endgroup

\section*{Declaration of AI assistance}
During the preparation of this work, the authors used OpenAI Codex to assist with the development of some mathematical proofs. The authors reviewed the AI-assisted material and take full responsibility for the correctness of the proofs and the content of the article.

\end{document}